\documentclass[12pt]{article}

\usepackage{amsthm,amsmath,amssymb,multirow,subfigure,makecell,booktabs,array,url,bm,mathtools,wrapfig,lipsum,mathrsfs,relsize,titling,dsfont,epstopdf,bbm,color,caption, comment,enumitem, geometry}
\usepackage[pdftex]{graphicx}
\usepackage{xcolor}
\usepackage{algorithm, algpseudocode}
\usepackage[export]{adjustbox}[2011/08/13]

\setlist{leftmargin=5mm}
\usepackage[colorlinks,linkcolor=black,citecolor=blue,urlcolor=blue]{hyperref}

\newtheoremstyle{mystyle}
{}
{}
{}
{}
{\bfseries}
{.}
{ }
{\thmname{#1}\thmnumber{ #2}\thmnote{ (#3)}}

\usepackage{natbib}
\theoremstyle{mystyle}
\newtheorem{theorem}{Theorem}[section] 

\newtheorem{lemma}[theorem]{Lemma} 
\newtheorem{corollary}[theorem]{Corollary}

\newtheorem{proposition}[theorem]{Proposition}
\newtheorem{assumption}[theorem]{Assumption}
\newtheorem{remark}[theorem]{Remark}

\usepackage{thm-restate, cleveref}
\usepackage[page,toc,titletoc,title]{appendix}

\numberwithin{equation}{section}

\newcommand{\Fcal}{\mathcal{F}}

\newcommand{\Ncal}{\mathcal{N}}

\newcommand{\Wcal}{\mathcal{W}}

\newcommand{\EE}{\mathbb{E}}

\newcommand{\PP}{\mathbb{P}}
\newcommand{\QQ}{\mathbb{Q}}
\newcommand{\RR}{\mathbb{R}}

\newcommand{\argmin}{\mathop{\arg\min}}

\newcommand{\dd}{\text{d}}

\newcommand{\diam}{\text{diam}}

\newcommand{\op}{\text{op}}

\newcommand{\var}{\text{Var}}

\providecommand{\inner}[1]{\langle#1\rangle}

\newcommand{\sign}[1]{\text{sign}(#1)}

\newcommand{\vertiii}[1]{{\left\vert\kern-0.25ex\left\vert\kern-0.25ex\left\vert #1 
    \right\vert\kern-0.25ex\right\vert\kern-0.25ex\right\vert}}

\makeatletter
\def\keywords{\xdef\@thefnmark{}\@footnotetext}
\makeatother

\begin{document}

\title{Statistical Convergence Rates of Optimal Transport Map Estimation between General Distributions}

\author{
{Yizhe Ding} 
\and
{Runze Li} 
\and
{Lingzhou Xue}
}
\date{Department of Statistics, The Pennsylvania State University}

\maketitle

\begin{abstract}
This paper studies the convergence rates of optimal transport (OT) map estimators, a topic of growing interest in statistics, machine learning, and various scientific fields. 
Despite recent advancements, existing results rely on regularity assumptions that are very restrictive in practice and much stricter than those in Brenier's Theorem, including the compactness and convexity of the probability support and the bi-Lipschitz property of the OT maps. We aim to broaden the scope of OT map estimation and fill this gap between theory and practice. Given the strong convexity assumption on Brenier's potential, we first establish the non-asymptotic convergence rates for the original plug-in estimator without requiring restrictive assumptions on probability measures. Additionally, we introduce a sieve plug-in estimator and establish its convergence rates without the strong convexity assumption on Brenier's potential, enabling the widely used cases such as the rank functions of normal or $t$-distributions. We also establish new Poincaré-type inequalities, which are proved given sufficient conditions on the local boundedness of the probability density and mild topological conditions of the support, and these new inequalities enable us to achieve faster convergence rates for Donsker function class. Moreover, we develop scalable algorithms to efficiently solve the OT map estimation using neural networks and present numerical experiments to demonstrate the effectiveness and robustness.
\end{abstract}


\section{Introduction}\label{intro}

Given two probability measures $P$ and $Q$ on $\RR^d$, the \textit{Monge problem}, first formalized by \cite{monge1781memoire}, seeks to find a transport map $T: \RR^d\to \RR^d$ that transports $P$ to $Q$ with the minimal transportation cost, which can be formulated as
\begin{equation}
    \min_{T:\RR^d\to\RR^d} \int_{\RR^d} \|x-T(x)\|_2^2\, P(\dd x),\qquad \text{subject to}\quad T_\# P=Q,
\end{equation}
where $T_\# P=Q$ indicates that $T$ is a \textit{push-forward} from $P$ to $Q$, meaning that for any Borel set $A\subseteq\RR^d$, $Q(A)=P(T^{-1}(A))$. If a map denoted by $T_0$ attains the minimum of the Monge problem, $T_0$ is called an \textit{optimal transport (OT) map} from $P$ to $Q$.
However, the Monge problem does not always admit a solution. Let $\Pi(P,Q)$ be the collection of probability measures on $\RR^d\times\RR^d$ with marginals $P$ and $Q$. The \textit{Kantorovich problem} \citep{kantorovich1942translocation, kantorovich2006translocation} relaxes the notion of push-forward and seeks to find a probability measure $\pi\in\Pi(P,Q)$ to solve
\begin{equation}
    \min_{\pi} \iint_{\RR^d\times\RR^d} \|x-y\|_2^2\, \pi(\dd x, \dd y),\qquad\text{subject to}\quad \pi\in\Pi(P,Q).
\end{equation}
The minimal value of the Kantorovich problem is known as the squared $2$-Wasserstein distance between $P$ and $Q$, denoted by $\Wcal_2^2(P,Q)$.
Brenier's theorem \citep{brenier1987decomposition, brenier1991polar}  bridges the Monge and Kantorovich problems by showing that the OT map $T_0$ can be recovered from the Kantorovich problem and represented as $T_0=\nabla\varphi_0$ almost everywhere, which is the gradient of a convex function $\varphi_0$ known as \textit{Brenier's potential}.

Over the past decade, optimal transport has emerged as an important concept \citep{villani2003topics, villani2009optimal, santambrogio2015optimal}, with notable applications across various scientific disciplines. For instance, optimal transport has shown great promise in advancing methods for single-cell data analysis \citep{schiebinger2019optimal,cao2022unified,bunne2023learning}. In statistics, optimal transport offers a novel framework for studying multivariate ranks and quantiles. In the univariate setting, the rank and quantile functions correspond exactly to the OT maps between the distributions of interest and the uniform measure on the unit interval $(0,1)$.  Multivariate ranks and quantiles based on optimal transport retain properties analogous to their univariate counterparts \citep{chernozhukov2017monge, hallin2021distribution, ghosal2022multivariate}. These developments have further inspired new approaches to distribution-free testing \citep{deb2023multivariate, huang2023multivariate} and copula-based graphical modeling \citep{zhang2024copula}.

These successful applications of optimal transport highlight the significance of studying the theoretical properties of optimal transport estimators. \cite{hutter2021minimax} introduced three OT map estimators and studied the minimax estimation of smooth OT maps, assuming that $P$ and $Q$ are compactly supported on hypercubes. Their framework required the OT map to be $\alpha$-Hölder smooth and Brenier's potential to be strongly convex. Building on this foundation, \cite{manole2021plugin} analyzed three OT map estimators introduced by \cite{hutter2021minimax}. While \cite{manole2021plugin} allowed $P$ to have a Hölder smooth density in their empirical estimator, they imposed Hölder smoothness on the density of $Q$ in their wavelet-based estimator, and their 1-NN estimator introduced additional topological constraints on the supports of $P$ and $Q$.
However, the assumptions of convexity and compactness on the probability support of $P$ and $Q$ in \cite{hutter2021minimax} and \cite{manole2021plugin} restrict the study of OT maps and do not apply to OT maps between probability measures with non-convex or unbounded supports. This gap impedes the development of robust inference tools and slows progress in addressing unresolved statistical challenges using optimal transport. Notably, Brenier's theorem \citep{brenier1987decomposition, brenier1991polar} only requires $P$ and $Q$ to have finite second-order moments, indicating that exploring approaches that relax these restrictive assumptions is natural and necessary. Such extensions would broaden the applicability of optimal transport to more general and practical scenarios.

\cite{deb2021rates} advanced the study of OT map estimation by requiring only $P$ to be sub-Weibull and $Q$ to be compactly supported, with the OT map assumed to be Lipschitz in their Theorem 2.2. However, this estimator suffers from the curse of dimensionality. To address this issue, they reintroduced assumptions on the smoothness of the probability densities of $P$ and $Q$, as well as the compactness and convexity of their supports, enabling the use of kernel density estimation in their Theorem 2.5. While this mitigates the dimensionality issue, it limits the generality of their results and introduces additional bias in practice due to the selection of kernel and bandwidth. Moreover, as a barycentric projection, their OT map estimator cannot be computed explicitly, and its polynomial computational complexity makes the approximated estimator inefficient for scaling with increasing sample sizes.

Following \cite{hutter2021minimax}, \cite{gunsilius2022convergence} explored constructing the OT map estimator based on the semi-dual form of the Kantorovich problem. \cite{gunsilius2022convergence} did not assume the smoothness or strong convexity of Brenier's potential. Instead, the analysis relied on the assumption that $P$ satisfies the Poincaré inequality, a sufficient condition for the local identification condition required to apply the empirical process theory. However, similar to prior work, the compactness and convexity of the supports of $P$ and $Q$ were still necessary to the analysis of \cite{gunsilius2022convergence}. Additionally, the boundedness and local Hölder smoothness of the probability densities of $P$ and $Q$ were required in \cite{gunsilius2022convergence}.

\cite{divol2022optimal} extended the use of Poincaré inequality to allow the supports of $P$ and $Q$ to be unbounded, provided they are restricted to hyperballs in $\RR^d$. However, in cases where $P$ does not have bounded support, they assumed that Brenier's potential $\varphi_0$ is $(\alpha, a)$-convex, which excludes scenarios where the OT map corresponds to the rank function of a normal distribution. while the Poincaré inequality is widely applicable in probability and functional analysis and holds for many common probability measures, it implies that $P$ must be sub-exponential \citep{bobkov1997poincare}. This excludes many distributions of practical interest, such as the $t$-distribution and Pareto distribution, which are not even sub-Weibull. Addressing these heavy-tailed distributions is essential for advancing optimal transport theory to more diverse and practical applications.

Computational considerations are also critical for advancing OT map estimation. To bridge the statistical-computational gap, \cite{muzellec2021near} leveraged the sum-of-squares (SoS) tight reformulation of OT introduced by \cite{vacher2021dimension}. This approach transformed the semi-dual form of the Kantorovich problem into an unconstrained convex program. However, their method required the strong assumption that both Brenier's potential $\varphi_0$ and its convex conjugate $\varphi_0^\ast$ belong to the Sobolev space $H^{m+2}(\RR^d)$ with $m>d+1$. Moreover, \cite{divol2022optimal} showed that solving the OT estimator using the semi-dual formulation of the Kantorovich problem with a general function class is NP-hard. These challenges highlight the pressing need for practical algorithmic advancements that reduce computational complexity while maintaining theoretical robustness and narrowing the statistical-computational gap.

\subsection*{Main Contributions}

In this work, we significantly relax the assumptions currently made in the literature to derive statistical convergence rates of OT map estimators based on the semi-dual formulation of the Kantorovich problem. We summarize the assumptions in existing literature alongside our contributions in Table~\ref{table: assumptions of OT estimation, measure P} and Table~\ref{table: assumptions of OT estimation, potential} for clarity. Our approach extends the applicability of OT map estimation from sub-exponential measures or distributions with compact supports to general sub-Weibull distributions and even distributions with only $(2+\epsilon)$-order moments, nearly closing the gap between existing methods and the theoretical limit established by Brenier's Theorem. Our main theoretical and methodological contributions are summarized below:

\begin{enumerate}
    \item \textbf{Non-Asymptotic Rates with Minimal Requirements on Probability Measures:} Under the $(\alpha,a)$-convexity assumption of Brenier's potential, Theorem \ref{thm: general function class with strongly convex, main text} in Section \ref{sec: General Function Class} establishes non-asymptotic convergence rates for OT map estimators when $\Fcal$ is a general function class. To the best of our knowledge, this is the first result in the literature to provide non-asymptotic convergence rates for OT map estimators with minimal requirements on probability measures. 
    
    Prior works, such as \cite{deb2021rates} and \cite{divol2022optimal}, addressed unbounded probability measures but had key limitations: \cite{deb2021rates} assumed $Q$ was bounded, and \cite{divol2022optimal} restricted the domain to hyperballs and required $P$ to be sub-exponential.  
    In contrast, our approach does not require $P$ and $Q$ to have compact or convex supports, nor do we impose smooth density assumptions. This significantly broadens the applicability of OT map estimation and provides robust theoretical guarantees for the convergence rates of OT-based multivariate ranks and quantiles. Existing works have been more restrictive: for example, \cite{chernozhukov2017monge} established only consistency for OT-based multivariate ranks and quantiles, while \cite{ghosal2022multivariate} derived convergence rates exclusively for probability measures with compact and convex supports.
    
    \item \textbf{Improved Asymptotic Rates for Donsker Function Classes and New Poincaré-Type Inequalities:} Poincaré inequalities have been instrumental in establishing convergence rates for OT map estimators in the existing literature. Theorem \ref{thm: donsker function class with strongly convex, main text} in Section \ref{sec: Donsker Function Class} extends the results of Theorem \ref{thm: general function class with strongly convex, main text}, focusing on cases where $\Fcal$ is a Donsker function class and $P$ satisfies certain Poincaré-type inequalities. While this result pertains to asymptotic convergence rates rather than non-asymptotic ones, it is notable for enabling faster, and even parametric, convergence rates. Such rates are achievable by common statistical models, including parametric models, nonlinear models using Reproducing Kernel Hilbert Spaces (RKHS) with kernels that have exponentially decaying spectra, and neural networks. 
    
    These improved asymptotic rates are driven by the use of Poincaré inequalities, as seen in prior works like \cite{gunsilius2022convergence} and \cite{divol2022optimal}. However, those works required $P$ to be sub-exponential as a necessary condition for applying Poincaré inequalities. In contrast, our results leverage our newly developed Poincaré-type inequalities in Section \ref{sec: Poincaré-type Inequalities}, which extend the applicability to sub-Weibull and polynomial-tailed distributions, significantly broadening the scope and practical utility of these advancements.
    
    Unlike traditional Poincaré inequalities, our newly developed Poincaré-type inequalities are not dependent on the tail thickness of the probability measures. Instead, they can be established based on the local boundedness of the probability density and mild topological properties of the probability support. These newly developed Poincaré-type inequalities may become versatile probabilistic tools with potential applications in non-parametric statistics and other areas of statistical research.

    \item \textbf{A New Sieve Estimator Eliminating the $(\alpha,a)$-Convexity Assumption:} Inspired by sieve estimates \citep{shen1994convergence}, we introduce a novel sieve plug-in OT map estimator that eliminates the need for the $(\alpha,a)$-convexity assumption of Brenier's potential required in Theorem \ref{thm: general function class with strongly convex, main text} and Theorem \ref{thm: donsker function class with strongly convex, main text}. This $(\alpha,a)$-convexity condition is impractical in many cases; for instance, the rank function of the normal distribution, as an OT map from the normal distribution to the uniform measure on $(0,1)$, fails to be $(\alpha,a)$-convex for any $a\geq 0$. 
    
    Our new estimator addresses this limitation and, to the best of our knowledge, is the first in the literature to enable OT map estimation between an unbounded measure and a bounded one, or between a polynomial-tailed distribution and a sub-Weibull distribution. The refined results corresponding to Theorem \ref{thm: general function class with strongly convex, main text} and Theorem \ref{thm: donsker function class with strongly convex, main text} are presented in Theorem \ref{thm: general function class without strongly convex, main text} and Theorem \ref{thm: donsker function class without strongly convex, main text}, respectively.
    
    \item \textbf{Approaching the Theoretical Limit of Brenier's Theorem:} While we have extended the convergence rates of plug-in OT map estimators from cases where $P$ has compact support to cases where $P$ has at least fourth-order moments in Sections \ref{sec: General Function Class}, \ref{sec: Donsker Function Class} and \ref{sec: Relaxing the Strong Convexity}, Brenier's Theorem only requires $P$ and $Q$ to have second-order moments. This gap arises because existing empirical process tools typically require functions to be $L^2(P)$-integrable \citep{koltchinskii2011oracle, boucheron2013concentration, vaart2023empirical}. 
    To bridge this gap, we develop a novel empirical process result for unbounded $L^1(P)$-integrable functions in Section \ref{sec: Beyond 4th-order Moments}, detailed in Lemma \ref{lemma: convergence of EP with L^1 integrable functions}. To the best of our knowledge, this is the first result in the literature specifically addressing unbounded $L^1(P)$-integrable functions. This advancement allows us to derive convergence rates for OT map estimators in Theorem \ref{thm: convergence rate of semi-discrete optimal transport, second moment, main text}, even when $P$ and $Q$ lack fourth-order moments, thereby narrowing the gap to the theoretical limit set by Brenier's Theorem.
\end{enumerate}

Additionally, we introduce scalable algorithms that leverage neural networks for efficiently solving the OT map estimation problem. Both the original plug-in estimator and the newly introduced sieve plug-in estimator are computationally efficient and scalable, with their effectiveness and robustness validated through numerical experiments. Notably, the sieve estimator excels in robustness to heavy-tailed distributions. Simulations demonstrate that it consistently achieves lower $L^2(P)$ losses compared to the original estimator, particularly in heavy-tailed settings such as $t$-distributions.

\subsection*{Organization}

The structure of this work is as follows. We begin by introducing useful notations in the remainder of Section \ref{intro}. The preliminaries on optimal transport are provided in Section \ref{Sec Background of Optimal Transport}. In Section \ref{sec: main results}, we present the convergence rates of two plug-in OT map estimators under various settings. The newly developed Poincaré-type inequalities are presented in Section \ref{sec: Poincaré-type Inequalities}. In Section \ref{sec: numerical experiments}, we review existing OT map estimators from a computational perspective, introduce our proposed algorithm, and demonstrate the numerical performance of our estimators. Section \ref{conclusion} includes a few concluding remarks. More details about the algorithms, numerical experiments, and visualization are presented in Appendix A, and proofs are presented in the supplementary materials. 

\subsection*{Notations}

Throughout this paper, the following notations and terminologies are used. Let $\RR^d$ be the $d$-dimensional Euclidean space, equipped with Euclidean norm $\|\cdot\|_2$ and inner product $\inner{\cdot,\cdot}$. For $x_0\in\RR^d$ and $r>0$, denote $B(x_0,r)=\{x\in\RR^d:\|x-x_0\|_2\leq r\}$ as the ball in $\RR^d$. Given a metric space $(X,d_X)$, 
the diameter of $D\subseteq X$ is defined as $\diam(D)=\sup\limits_{x,y\in D}d(x,y)$. $D\subset\RR^d$ is called a Lipschitz domain with constants $r,L>0$, if for any boundary point $x_0\in\partial D$, after a transformation of coordinates if necessary, there is a Lipschitz function $f:\RR^{d-1}\to\RR$ satisfying $\|f\|_{L^\infty(B(x_0,r))}+r\|f\|_\text{Lip}\leq Lr$ and $D\cap B(x_0,r)=\{x=(x_1,x_2,\cdots,x_d)\in B(x_0,r): x_d>f(x_1,x_2,\cdots,x_{d-1})\}$. 

We say $a\lesssim b$, if there exists a constant $C$ independent of $a$ and $b$, such that $a\leq C b$. We say $a\asymp b$ if both $a\lesssim b$ and $b\lesssim a$ hold. We say $a\leq_{\log n}b$ if there exists a function $C(\log n)$ depending only on $\log n$ such that $a\leq C(\log n)\cdot b$. Notations $\lesssim_{\log n}$ and $\lesssim_{\log\log}$ are defined correspondingly. $a\lor b$ and $a\land b$ represent the maximal and minimal values between $a$ and $b$, respectively. For $x\in \RR$, denote $\log_+(x)=1\lor \log(x)$. For $x\in\RR^d$, denote $\inner{x}=(1+\|x\|_2)$. 

Given a probability measure $P$, the expectation of function $f$ is $\EE_P[f(X)]=Pf=\int f(x)\, P(\dd x)$. Besides, we use $\EE^\ast$ to represent outer expectation. 
We assume that $X_1,\cdots,X_n$ are i.i.d. copies from $P$, and $Y_1,\cdots,Y_N$ are i.i.d. copies from $Q$. The corresponding empirical measures are denoted as $\PP_n=\frac{1}{n}\sum_{i=1}^n \delta_{X_i}$ and $\QQ_N=\frac{1}{N}\sum_{j=1}^N \delta_{Y_j}$ respectively, where $\delta_\cdot$ is the Dirac measure. 
A random vector $X$ is called sub-Weibull with parameters $(\theta,K)$ (or just $\theta$ for simplicity), if for any $t>0$, $\PP(\|X\|_2\geq t)\leq 2\cdot \exp(-t^\theta/K^\theta)$. 

Let $C^m(\RR^d)$ be the space of $m$-times differentiable functions defined on $\RR^d$. The gradient and Hessian of $f\in C^2(\RR^d)$ are denoted as $\nabla f$ and $\nabla^2 f$, respectively. Denote $W^{1,2}(P)=H^1(P)$ as the $(1,2)$-Sobolev space with norm $\|f\|_{W^{1,2}(P)}=\|f\|_{L^2(P)}+\|\nabla f\|_{L^2(P)}$.
A function $\varphi\in C^2(\RR^d)$ is said to be $(\beta,b)$-smooth for some $\beta>0$ and $b\geq 0$, if $\|\nabla^2 \varphi(x)\|_\op\leq \beta \inner{x}^b$ for all $x\in\RR^d$, where the operator norm $\|\cdot\|_\op$ represents the largest absolute values of eigenvalues. $\varphi$ is said to be $(\alpha,a)$-convex for some $\alpha>0$ and $a\geq 0$, if $\lambda_{\min{}}{(\nabla^2 \varphi(x))}\geq \alpha \inner{x}^a$ for all $x\in\RR^d$, where $\lambda_{\min{}}(A)$ is the smallest eigenvalue of matrix $A$. 

\section{Preliminaries}\label{Sec Background of Optimal Transport}

\subsection{Monge and Kantorovich Problems}\label{Sec Monge and Kantorovich Problems}

Given two probability measures $P$ and $Q$ on $\RR^d$, a map $T:\RR^d\to \RR^d$ is called a push-forward from $P$ to $Q$, denoted as $T_\# P=Q$, if for any Borel measurable set $A\subseteq\RR^d$, its pre-image of $T$ satisfies $Q(A)=P(T^{-1}(A))$. It is possible that there is no push-forward between two probability measures. For instance, a discrete measure cannot be transported to an absolutely continuous probability measure, such as a normal distribution, because a push-forward, as a function, cannot be one-to-many. 

Even if push-forwards from $P$ to $Q$ exist, they may not be unique. For example, both $T_1(x) = x + 1$ and $T_2(x) = -x + 1$ are push-forwards from normal distribution $P=N(0,1)$ to $Q=N(1,1)$. Therefore, it is natural to ask, among all possible push-forwards from $P$ to $Q$, which one is ``optimal'' in the sense of minimizing the the transportation cost from $P$ to $Q$. This \textit{Monge problem} can be formulated as the following minimization program:      
\begin{equation}\label{eq: monge problem definition}
    \min_{T:\RR^d\to\RR^d} \int_{\RR^d} \|x-T(x)\|_2^2\, P(\dd x),\qquad \text{subject to }\quad T_\# P=Q,
\end{equation}
where $\|x-T(x)\|_2^2$ represents the transportation cost between $x$ and $T(x)$. Its solution, $T_0$, is referred to as the \textit{optimal transport map} (OT map) from $P$ to $Q$.

Regardless of the infeasibility issue of the Monge problem, the constraint $T_\# P=Q$ is also difficult to handle. These problems remained unresolved until the relaxation proposed by \cite{kantorovich1942translocation, kantorovich2006translocation}. Instead of optimizing over the collection of push-forwards from $P$ to $Q$, Kantorovich sought to find a product measure on $\RR^d\times\RR^d$, with marginals $P$ and $Q$, that minimizes the transportation cost between its two marginal measures. This \textit{Kantorovich problem} can be expressed as 
\begin{equation}\label{eq: kantorovich problem definition}
    \min_{\pi} \iint \|x-y\|_2^2\, \pi(\dd x, \dd y),\qquad\text{subject to }\quad \pi\in\Pi(P,Q).
\end{equation}
Here, $\Pi(P,Q)$ represents the collection of couplings with $P$ and $Q$ as the marginals. 

By considering couplings of $P$ and $Q$, Kantorovich problem addresses the infeasibility issue in the Monge problem, as push-forwards can represent special couplings between $P$ and $Q$. In fact, if $T_\# P=Q$, then $(\text{id}\times T)_\# P$ is a coupling of $P$ and $Q$. Moreover, since the collection $\Pi(P,Q)$ is convex, and the objective function is linear with respect to $\pi$, the Kantorovich problem is a linear program, albeit an infinite-dimensional one. 

Lastly, we derive the semi-dual form of Kantorovich problem from the original one in Equation \eqref{eq: kantorovich problem definition}. Denote $\inner{\cdot,\cdot}$ as the inner product in $\RR^d$. By substituting $\|x-y\|_2^2=\|x\|_2^2+\|y\|_2^2-2\inner{x,y}$ into the Kantorovich problem in Equation \eqref{eq: kantorovich problem definition}, we get $$\min_{\pi\in\Pi(P,Q)}\iint \|x-y\|_2^2 \pi(\dd x, \dd y)=\EE_P[\|X\|_2^2] +\EE_Q[\|Y\|_2^2]-2\max_{\pi\in\Pi(P,Q)}\iint \inner{x,y}\pi(\dd x, \dd y).$$ As $\EE_P[\|X\|_2^2] +\EE_Q[\|Y\|_2^2]$ is independent of the coupling $\pi$, we can focus on the third term $\max_{\pi}\iint \inner{x,y}\pi(\dd x, \dd y)$. Since it is linear with respect to $\pi$, its dual form reads 
\begin{equation}\label{eq, semi-dual of Kantorovich definition}
    \min_{\varphi\in L^1(P)} P\varphi+Q\varphi^\ast,\quad\text{where}\quad \varphi^\ast(y)=\sup_{x\in\RR^d}\inner{x,y}-\varphi(x).
\end{equation}
$\varphi^\ast$ is the convex conjugate (a.k.a. Legendre–Fenchel transformation) of $\varphi$ and Equation \eqref{eq, semi-dual of Kantorovich definition} is called the \textit{semi-dual form of the Kantorovich problem}.

\subsection{Brenier's Theorem}\label{sec: Brenier's Theorem}

The Kantorovich problem can be regarded as a relaxation of the Monge problem. So, one may naturally ask if they are equivalent and whether we can recover the optimal transport from the Kantorovich problem. These questions were addressed in \cite{brenier1987decomposition, brenier1991polar}:
\begin{theorem}[Brenier's Theorem \citep{brenier1987decomposition, brenier1991polar}]\label{thm: Brenier theorem}
    Assume $P$ and $Q$ are two Borel probability measures on $\RR^d$ such that $P$ is absolutely continuous with respect to the Lebesgue measure, and $\EE_P\|X\|_2^2<\infty$, $\EE_Q\|Y\|_2^2<\infty$. Then there exists a convex function $\varphi_0:\RR^d\to\RR\cup\{+\infty\}$ such that its gradient $\nabla\varphi_0:\RR^d\to\RR^d$ is the unique push-forward from $P$ to $Q$ that arises as the gradient of a convex function. Moreover, $\nabla\varphi_0$ uniquely minimizes the Monge problem in Equation \eqref{eq: monge problem definition}, $(\text{id}\times \nabla\varphi_0)_\# P$ uniquely minimizes the Kantorovich problem in Equation \eqref{eq: kantorovich problem definition} and $\varphi_0$ is the unique minimizer (up to adding a constant) of the semi-dual form of Kantorovich problem in Equation \eqref{eq, semi-dual of Kantorovich definition}. 
\end{theorem}

Brenier's Theorem is also known as the Polar Factorization Theorem, which states that each push-forward $T_\#P=Q$ can be uniquely decomposed as $T = \nabla\varphi_0 \circ s$, where $\nabla\varphi_0$ is the OT map from $P$ to $Q$, and $s$ is a measure-preserving map of $P$, i.e., for any Borel set $A$, $P(A) = P(s^{-1}(A))$. Hence, to study the transport from $P$ to $Q$, it is impossible for statistician to establish an identifiable model based solely on push-forwards, due to the existence of the measure-preserving map. 

Now, consider the example of $P=N(0,1)$ and $Q=N(1,1)$ mentioned earlier. Clearly, $T_1(x) = x + 1$, as the gradient of the convex function $\frac{1}{2}x^2 + x$, is the OT map from $P$ to $Q$. Meanwhile, the polar decomposition of the push-forward $T_2 = -x + 1$ can be written as $T_2 = T_1 \circ s$, where $s(x) = -x$ is a measure-preserving map of $P$. In the following, we shall refer to $\varphi_0$ as the Brenier's potential, and $\nabla\varphi_0$ as the OT maps.

Optimal transport can be viewed as a natural generalization of cumulative distribution functions (CDFs) or rank functions in $\RR$. Consider a univariate absolutely continuous random variable $X\sim P$ with CDF $\mathbf{F}(\cdot)$. Remarkably, $\mathbf{F}(X)$ follows $U(0,1)$, the uniform measure on interval $(0,1)$. Since $\mathbf{F}(\cdot)$, as a non-decreasing function, can be represented as the derivative of some convex function, Brenier's Theorem implies that $\mathbf{F}(\cdot)$ is the OT map from $P$ to $U(0,1)$. Similarly, the quantile function, $\mathbf{Q}(\cdot)=\mathbf{F}^{-1}(\cdot)$, is the OT map from $U(0,1)$ to $P$.

\section{Main Results}\label{sec: main results}

In this section, we investigate the convergence rates of plug-in OT map estimators across various settings. We begin with an overview of two plug-in estimators and discuss the assumptions made on probability measures in the literature in Section \ref{sec: overview}. Sections \ref{sec: General Function Class} and \ref{sec: Donsker Function Class} present the convergence rates for the first estimator under the $(\alpha,a)$-convexity assumption. In Section \ref{sec: L_infty Rates with Bounded Supports}, we establish a connection between semi-local uniform convergence rates and $L^2(P)$ rates, facilitating the derivation of $L^\infty$ convergence rates under bounded support condition. The convergence rates of the second estimator, which does not rely on the $(\alpha,a)$-convexity assumption, are introduced in Section \ref{sec: Relaxing the Strong Convexity}. Finally, in Section \ref{sec: Beyond 4th-order Moments}, we examine both estimators when the probability measure $P$ may lack fourth-order moments.

\subsection{Overview of Plug-in OT Map Estimators}\label{sec: overview}

Brenier's Theorem provides a foundational framework for estimating the OT map using a plug-in approach through the semi-dual form of the Kantorovich problem as defined in Equation \eqref{eq, semi-dual of Kantorovich definition}. Specifically, given two empirical measures $\PP_n$ and $\QQ_N$ from distributions $P$ and $Q$ respectively, which may not be independent, the first plug-in estimator of the Brenier's potential is formulated as
\begin{equation}\label{eq: discrete-discrete OT}
    \hat{\varphi}_{n,N} \in \argmin_{\varphi\in\Fcal}\PP_n\varphi+\QQ_N \varphi^\ast,\quad\text{where}\quad \varphi^\ast(y)=\sup_{x\in\RR^d} \inner{x,y}-\varphi(x).
\end{equation}
Consequently, the OT map estimator from $P$ to $Q$ is defined as $\nabla \hat{\varphi}_{n,N}$. As both $P$ and $Q$ are unknown in Equation \eqref{eq: discrete-discrete OT} and are replaced with their empirical counterparts, $\nabla \hat{\varphi}_{n,N}$ is termed the \textit{discrete-discrete estimator}. 

While prior studies, such as \cite{gunsilius2022convergence} and \cite{divol2022optimal}, have primarily studied this estimator, strong convexity condition (or $(\alpha,a)$-convexity) of Brenier's potential is needed to establish its convergence rates when $P$ may have a unbounded support. To address scenarios where the $(\alpha,a)$-convexity does not hold, we introduce an alternative ``sieve'' plug-in estimator inspired by \cite{shen1994convergence}:
\begin{equation}\label{eq: new discrete-discrete estimator, overview}
    \tilde\varphi_{n,N}\in\argmin_{\varphi\in\Fcal}\PP_n\varphi+\QQ_N\varphi^{\ast,(n)},\quad\text{where}\quad \varphi^{\ast,(n)}(y)=\sup_{x\in B(0,M_n)}\inner{x,y}-\varphi(x).
\end{equation}
In this formulation, we restrict the supremum's search domain from $\RR^d$ to the bounded hyperball $B(0,M_n)$. This modification effectively serves as a pseudo-support for $P$, allowing us to bypass the necessity of the $(\alpha,a)$-convexity assumption. The motivation behind it shall be further discussed in Section \ref{sec: Relaxing the Strong Convexity}. For both theoretical analysis and practical implementations, we set $M_n=\max_{1\leq i\leq n}\|X_i\|_2$, representing the maximum norm among all samples in $\PP_n$. We hypothesize that alternative choices of $\{M_n\}_n$ or different shapes for the supremum's search domain could potentially yield better convergence rates, and we leave this exploration to interested readers.

It is noteworthy that the working function class $\Fcal$ does not need to consist solely of convex functions, despite Brenier's potential being convex. This flexibility aligns with approaches in \cite{divol2022optimal} and allows practitioners to leverage a wide range of models, including wavelets, Reproducing Kernel Hilbert Spaces (RKHS), and neural networks, to estimate the OT maps.

Given that the true Brenier's potential $\varphi_0$ may not belong to the working function class $\Fcal$, we define $\overline\varphi\in\argmin_{\varphi\in\Fcal}P\varphi+Q\varphi^\ast$ as the projection of $\varphi_0$ within $\Fcal$. Then $\|\nabla\hat\varphi_{n,N}-\nabla\varphi_0\|_{L^2(P)}$ can be decomposed into the approximation error and the estimation error:
$$\|\nabla\hat\varphi_{n,N}-\nabla\varphi_0\|_{L^2(P)}\leq \|\nabla\overline\varphi-\nabla\varphi_0\|_{L^2(P)}+\|\nabla\hat\varphi_{n,N}-\nabla\overline\varphi\|_{L^2(P)}.$$ Since the approximating error $\|\nabla\overline\varphi-\nabla\varphi_0\|_{L^2(P)}$ depends on the relation between $\Fcal$ and $\varphi_0$, our main focus is on analyzing the estimation error $\|\nabla\hat\varphi_{n,N}-\nabla\overline\varphi\|_{L^2(P)}$.

Alternatively, it is referred to as a \textit{semi-discrete estimator} when either $P$ or $Q$ is known:
\begin{equation}\label{eq: semi-discrete OT}
    \hat{\varphi}_{n}\in\argmin_{\varphi\in\Fcal}\PP_n\varphi+Q\varphi^\ast,\qquad \hat{\varphi}_{N}\in\argmin_{\varphi\in\Fcal}P\varphi+\QQ_N \varphi^\ast.
\end{equation}

We will investigate the convergence rates of the discrete-discrete estimators from Equation \eqref{eq: discrete-discrete OT} and Equation \eqref{eq: new discrete-discrete estimator, overview}. Based on these results, one can get the convergence rates for semi-discrete estimators by tracing the terms involving $n$ and $N$ in the convergence rates. 

\begin{table}[ht]
    \centering
    \begin{tabular}{c|cc|cc}
        \hline
        & \multicolumn{2}{c|}{Supports} & \multicolumn{2}{c}{Density} \\ \cline{1-5}
        & Boundedness & Shape & Boundedness & Smoothness \\ \hline
        \cite{hutter2021minimax} & $P,Q$ & Hypercube & $P$ & ---  \\
        \cite{manole2021plugin} & $P,Q$ & Hypercube & $Q$ & $Q$ \\
        Theorem 2.2 in \cite{deb2021rates} & $Q$ & --- & --- & --- \\
        Theorem 2.5 in \cite{deb2021rates} & $P,Q$ & Convex & $P,Q$ & $P,Q$ \\
        \cite{muzellec2021near} & $P,Q$ & Convex & $P,Q$ & $P,Q$ \\
        \cite{gunsilius2022convergence} & $P,Q$ & Convex & $P,Q$ & $P,Q$ \\
        \cite{divol2022optimal} & --- & Hyperball & --- & ---\\
        Our Results & --- & --- & --- & --- \\ \hline
    \end{tabular}
    \caption{Assumptions on the probability measures in the literature. ``Hypercube'' means $[0,1]^d$, while ``hyperball'' refers to $B(0,R)$ with $R\in[0,\infty]$. ``Bounded density'' indicates that the density function is bounded away from $0$ and $\infty$, while ``smooth density'' means the probability density is Hölder smooth (locally Hölder smooth for \cite{gunsilius2022convergence}, Sobolev smooth for Theorem 2.5 in \cite{deb2021rates}).  
    In \cite{manole2021plugin}, Corollary 8 does not require the smoothness of probability densities of $P$ and $Q$, but it is subject to the curse of dimensionality; the density of $Q$ is required to be $(\alpha-1)$-Hölder smooth for their wavelet estimator in Theorem 10; additional topological assumptions on the supports of $P$ and $Q$ are required for their 1NN estimator in Proposition 15. In Theorem 2.2 of \cite{deb2021rates}, they assume $P$ is sub-Weibull, $Q$ is compactly supported, and the OT map is strongly convex. 
    For our work, Theorem \ref{thm: general function class with strongly convex, main text} and Theorem \ref{thm: general function class without strongly convex, main text} do not require any such assumptions for probability measures $P$ and $Q$ listed above. To achieve faster convergence rates, Poincaré-type inequalities are needed in Theorems \ref{thm: donsker function class with strongly convex, main text}, \ref{thm: donsker function class without strongly convex, main text} and \ref{thm: convergence rate of semi-discrete optimal transport, second moment, main text}.} 
    \label{table: assumptions of OT estimation, measure P}
\end{table}

Before presenting our results, we discuss the assumptions made in previous studies that we have successfully avoided, highlighting the innovations that distinguish our contributions. \cite{hutter2021minimax} laid the foundation for studying the convergence rates of OT maps estimators. In this seminal study, the authors investigated the discrete-discrete scenario with $n=N$, assuming that $P$ has a density supported on the hypercube $[0,1]^d$. They considered the OT map $\nabla\varphi_0$ as an $\alpha$-Hölder, Lipschitz continuous function, and assumed $\varphi_0$ to be strongly convex. 
    
Following \cite{hutter2021minimax}, \cite{manole2021plugin, muzellec2021near, deb2021rates, gunsilius2022convergence} extended this trajectory. Notably, almost all of these studies assumed the supports of $P$ and $Q$ to be compact and convex. The probability density function of $P$ was assumed to be bounded away from $0$ and $\infty$, and to be a smooth function with some degree of smoothness. For clarity and to avoid redundancy, we summarize the assumptions made in each study into Table~\ref{table: assumptions of OT estimation, measure P}.

Despite the convexity of the supports of $P$ and $Q$ is indeed a crucial condition for the existence of a continuous OT map (as discussed in Remark 3.8 of \cite{ghosal2022multivariate}), it is not strictly necessary. Moreover, there is a growing interest in investigating OT maps between probability measures with non-convex domains. For example, in the multivariate ranks and quantiles estimation, it is noteworthy that the Tukey halfspace depth-based estimator fails to deal with distributions with non-convex contours, such as banana-shaped distributions. In contrast, OT-based estimators demonstrate effectiveness in capturing such complex structures. We refer readers to Section 2.2 of \cite{chernozhukov2017monge} for further insights. 

In the reminder of this section, we explore the $L^2(P)$ convergence rates of $\nabla\hat{\varphi}_{n,N}$ from Equation \eqref{eq: discrete-discrete OT} across different contexts. In Section \ref{sec: General Function Class}, we consider the case when $P$ is a sub-Weibull or polynomial-tailed distribution, with $\Fcal$ being a general function class, as outlined in Theorem \ref{thm: general function class with strongly convex, main text}. In Section \ref{sec: Donsker Function Class}, we present improved convergence rates in Theorem \ref{thm: donsker function class with strongly convex, main text} for scenarios where $P$ additionally satisfies our Poincaré-type inequality and $\Fcal$ is a Donsker class. When $P$ is sub-Weibull, our convergence rates match those in \cite{divol2022optimal} with fewer assumptions. In Section \ref{sec: L_infty Rates with Bounded Supports}, we establish a connection between semi-local uniform convergence rates and $L^2(P)$ rates. Based on that, we provide $L^\infty$ convergence rates when $P$ has a bounded support. Then, in Section \ref{sec: Relaxing the Strong Convexity}, we discussed the sieve estimator, $\nabla\tilde\varphi_{n,N}$, designed to relax the $(\alpha,a)$-convexity assumption underlying the previous results. Its convergence rates are studied in Theorem \ref{thm: general function class without strongly convex, main text} and Theorem \ref{thm: donsker function class without strongly convex, main text} when $\Fcal$ is a general function class and a Donsker class, respectively. These results make us the first to study the convergence rates of the estimated rank function of a normal or $t$-distribution using optimal transport theory. Finally, in Section \ref{sec: Beyond 4th-order Moments}, we examine situations where $P$ and $Q$ lack fourth-order moments in Theorem \ref{thm: convergence rate of semi-discrete optimal transport, second moment, main text}, almost bridging the gap between the existing literature and the Brenier's Theorem. 

\subsection{General Function Class}\label{sec: General Function Class}

Let's introduce the following regularity assumptions:
\begin{assumption}[Existence of Brenier's potential]\label{assumption: existence of Brenier's potential, main}
    For probability measures $P$ and $Q$ on $\RR^d$,
    $\varphi_0\in C^2(\RR^d)$ is the Brenier's potential from $P$ to $Q$ with $\varphi_0(0)=0$.  
\end{assumption}

\begin{assumption}[Convexity and smoothness of potentials]\label{assumption: Convexity and Smoothness of Brenier's potential, main}
    $\varphi_0$ is the $(\beta,a)$-smooth, $(\alpha,a)$-convex Brenier's potential from $P$ to $Q$ for some $\alpha,\beta>0$ and a common $a\geq 0$. The working function class $\Fcal$ is made up of $(\beta,a)$-smooth functions with $\varphi(0)=0$ for all $\varphi\in\Fcal$.  
\end{assumption}

\begin{assumption}[Function class]\label{assumption: function class Fcal, main}
    Given a function class $\Fcal$ made up of $C^2(\RR^d)$ functions, we assume $\varphi(0)=0$ for all $\varphi\in\Fcal$ and $\{\|\nabla\varphi(0)\|_2: \varphi\in\Fcal\}$ is bounded. Besides, we assume there exist some $\gamma\geq 0$, $\gamma^\prime\geq 0$ and $D_\Fcal>0$ such that for every $h>0$, either the covering entropy, or the bracketing entropy of function space $\Fcal$ satisfies
    \begin{equation}\label{eq: L-2 covering entropy of Fcal}
        \EE[\log \Ncal(h, \Fcal,L^2(\PP_n))]\leq D_\Fcal\cdot h^{-\gamma}\cdot\log_+(1/h)^{\gamma^\prime}, \tag{{i}}
    \end{equation} 
    \begin{equation}\label{eq: L-2 bracketing entropy of Fcal}
        \text{or}\quad \log \Ncal_{[\ ]}(h, \Fcal,L^2(P))\leq D_\Fcal\cdot h^{-\gamma}\cdot\log_+(1/h)^{\gamma^\prime}. \tag{{ii}}
    \end{equation}
\end{assumption}

\begin{remark}[Comments to Assumption \ref{assumption: function class Fcal, main}]
    It is widely recognized that OT map estimation is notoriously affected by the curse of dimensionality (CoD). To mitigate this, previous studies, such as \cite{hutter2021minimax, manole2021plugin, deb2021rates, gunsilius2022convergence}, have assumed the smoothness of the OT maps or the density functions. However, when estimating the unknown density with kernel density estimation, selecting the kernel and bandwidth introduces additional biases in estimating Brenier's potential.

    In response to these challenges, our approach aligns with that of \cite{divol2022optimal}, which mitigates the CoD by constraining the complexity of the function class $\Fcal$. This strategy avoids imposing extra assumptions on $P$ and $Q$. While \cite{divol2022optimal} utilizes $L^\infty(\inner{\cdot}^{-\eta})$ covering entropy to accommodate unbounded probability supports, we considered $L^2$ counterparts for generality. 
    
    The covering/bracket entropy of the function class $\Fcal$ in Assumption \ref{assumption: function class Fcal, main} is dominated by two constants, $D_\Fcal$ and $\gamma$. The constant $D_\Fcal>0$ serves as an ``effective dimension'' of $\Fcal$, influencing both the approximation and estimation errors. Accordingly, we denote $\tilde{n}=n/D_\Fcal$ and $\tilde{N}=N/D_\Fcal$ as the ``effective sample sizes''. The constant $\gamma$ is pivotal in determining the complexity of $\mathcal{F}$; a larger $\gamma$ indicates a more complex $\Fcal$. We refer to $\Fcal$ as a general function class when $\gamma\geq 2$, and Donsker class when $\gamma<2$.  

    The introduction of $L^2$ covering/bracketing entropies provides a unified framework for deriving the convergence rates of the estimated OT maps across various function classes. In fact, even $L^\infty$ covering entropy is available for common non-parametric models, including RKHS (Lemma D.2 in \cite{yang2020function}) and deep neural networks (Lemma 5 in \cite{schmidt2020nonparametric}). These function classes may be dense in larger, more general function spaces like the Hölder, Sobolev and Besov spaces, while maintaining relatively lower model complexities. We refer readers to Section 2.7 of \cite{vaart2023empirical} for more examples on the covering entropies of different function spaces.
\end{remark}

Now, we are ready to present the convergence rate when $\Fcal$ is a general function space:
\begin{theorem}\label{thm: general function class with strongly convex, main text}
    Suppose Assumptions \ref{assumption: existence of Brenier's potential, main}, \ref{assumption: Convexity and Smoothness of Brenier's potential, main} and \ref{assumption: function class Fcal, main} \eqref{eq: L-2 covering entropy of Fcal} hold. $\Fcal$ is a general function class in the sense that $\gamma\geq 2$. Let $\overline\varphi$ be an arbitrary $(\alpha,a)$-convex potential in $\Fcal$. 
    
    \textbf{(i).} Suppose $P$ is a sub-Weibull probability measure with parameters $(\theta,K)$ for some $\theta,K>0$. Then, with probability $1-\delta_1-\delta_2$, 
    \begin{equation}
        \begin{aligned}
            &\|\nabla\hat\varphi_{n,N}-\nabla\varphi_0\|_{L^2(P)}^2
            \lesssim_{\log\log} \Bigl(S(\overline\varphi)-S(\varphi_0)\Bigr)\cdot\log_+\Bigl(\frac{1}{S(\overline\varphi)-S(\varphi_0)}\Bigr)^\frac{a(a+1)}{\theta}\\
            &\qquad+\tilde{n}^{-\frac{1}{\gamma}}\cdot\Bigl(\log\frac{\tilde{n}}{\delta_1/4}\Bigr)^{a+2+\frac{\gamma^\prime\lor 1}{2}+\frac{a(a+1)}{\theta}}
            +\tilde{N}^{-\frac{1}{\gamma}}\cdot\Bigl(\log\frac{\tilde{N}}{\delta_2/4}\Bigr)^{\frac{a+2}{a+1}+\frac{\gamma^\prime\lor 1}{2}+\frac{a(a+1)}{\theta}},
        \end{aligned}
    \end{equation}
    where $S(\varphi)=P\varphi+Q\varphi^\ast$.

    \textbf{(ii).} Suppose $P$ is a polynomial-tailed distribution in the sense that $\EE_P\|X\|_2^m<\infty$ with $m>(1+k)(2a+4)$ for some $k>0$. Then, with probability $1-n^{-k}-N^{-k}$, we have 
    \begin{equation}
        \begin{aligned}
            &\|\nabla\hat\varphi_{n,N}-\nabla\varphi_0\|_{L^2(P)}^2
            \lesssim\Bigl(S(\overline\varphi)-S(\varphi_0)\Bigr)^{\frac{m-2(a+1)}{m+(a+1)(a-2)}}\\
            &\qquad+(\tilde{n}\land \tilde{N})^{-(\frac{1}{\gamma}\land\frac{m-(1+k)(2a+4)}{2m})\frac{m-2(a+1)}{m+(a+1)(a-2)}}\cdot(\log (\tilde{n}\land \tilde{N}))^{\frac{\gamma^\prime+2}{2}\cdot\frac{m-2(a+1)}{m+(a+1)(a-2)}}.
        \end{aligned}
    \end{equation}
\end{theorem}

\begin{remark}[Comments to Theorem \ref{thm: general function class with strongly convex, main text}]
    Ignoring the logarithmic terms, Theorem \ref{thm: general function class with strongly convex, main text} suggests that the nonasymptotic estimation error for a sub-Weibull $P$ converges at the rate of $\tilde{n}^{-1/\gamma}+\tilde{N}^{-1/\gamma}$ when $\gamma\geq 2$, which aligns with the results presented in \cite{divol2022optimal}. 
    
    Compared to existing methods focusing on sub-Weibull measures or those with bounded support, we have successfully derived the convergence rates of the plug-in OT map estimator when $P$ has a polynomial-tail, addressing a significant gap in reaching the theoretical limit of $P$ presented by Brenier's Theorem (Theorem \ref{thm: Brenier theorem}). 

    Specifically, as the moment parameter $m$ increases, the convergence rate in the polynomial-tailed case approaches $\tilde{n}^{-1/\gamma}+\tilde{N}^{-1/\gamma}$. This demonstrates that having higher-order moments leads to improved estimation accuracy. Notably, many distributions of practical interest, such as the $t$-distribution and Pareto distribution, do not satisfy the sub-Weibull condition. However, Theorem \ref{thm: general function class with strongly convex, main text} accommodates these distributions provided that $\EE_P\|X\|_2^m<\infty$ with $m>2a+4$. This condition ensures that the $(\beta,a)$-smooth potentials $\varphi$'s are $L^2(P)$ integrable, thereby allowing the application of existing empirical process tools.
\end{remark}

\subsection{Donsker Function Class}\label{sec: Donsker Function Class}

In this section, we demonstrate that if the measure $P$ additionally satisfies the Poincaré-type inequality and $\Fcal$ is a Donsker class, then superior convergence rates compared to those in Theorem \ref{thm: general function class with strongly convex, main text} can be achieved.

\begin{assumption}[Poincaré-type inequality]\label{assumption: Poincare-type inequalities, main text}
    For any $\varphi_1,\varphi_2\in\Fcal$, 
    
    \textbf{(i).} when $P$ is sub-Weibull with parameter $\theta$, $$\var_P(\varphi_1-\varphi_2)\lesssim \|\nabla\varphi_1-\nabla\varphi_2\|_{L^2(P)}^2\cdot \log_+\Bigl(\frac{1}{\|\nabla\varphi_1-\nabla\varphi_2\|_{L^2(P)}}\Bigr)^{2/\theta}.$$

    \textbf{(ii).} When $P$ is a polynomial-tailed distribution with $\EE_P\|X\|_2^{2a+4+c}<\infty$ for some $c>0$, $$\var_P(\varphi_1-\varphi_2)\lesssim \|\nabla\varphi_1-\nabla\varphi_2\|_{L^2(P)}^2\lor\|\nabla\varphi_1-\nabla\varphi_2\|_{L^2(P)}^\frac{2c}{c+2}.$$
\end{assumption}

Our Poincaré-type inequalities closely resemble the global Poincaré inequality in Equation \eqref{eq: poincare inequality}. However, the key innovation lies in restricting the Poincaré-type inequalities to hold not for all Sobolev functions in $H^{1}(P)$, but only for differences $\varphi_1-\varphi_2$'s, where $\varphi_1$ and $\varphi_2$ are arbitrary potentials within the working function class $\Fcal$. 

One might question the generality and reasonableness behind this approach. We shall discuss the motivation behind our Poincaré-type inequalities and demonstrate their applicability to measure $P$ under several mild conditions in Section \ref{sec: Poincaré-type Inequalities}. Specifically, we will show that these inequalities hold for the arbitrary differences of $(\beta,b)$-smooth functions $\varphi_1-\varphi_2$, if the density of $P$ is locally bounded and its support satisfies certain topological conditions. 

\begin{theorem}\label{thm: donsker function class with strongly convex, main text}
    Suppose Assumptions \ref{assumption: existence of Brenier's potential, main}, \ref{assumption: Convexity and Smoothness of Brenier's potential, main}, \ref{assumption: function class Fcal, main} \eqref{eq: L-2 bracketing entropy of Fcal} and \ref{assumption: Poincare-type inequalities, main text} hold. $\Fcal$ is a Donsker function class in the sense that $0\leq\gamma<2$. Then, for any $(\alpha,a)$-convex $\overline\varphi\in\Fcal$,

    \textbf{(i).} when $P$ is sub-Weibull with parameters $(\theta,K)$ for some $\theta,K>0$. 
    \begin{equation}
        \begin{aligned}
            &\EE\|\nabla\hat\varphi_{n,N}-\nabla\varphi_0\|_{L^2(P)}^2
            \lesssim_{\log\log} \Bigl(S(\overline\varphi)-S(\varphi_0)\Bigr)\cdot\log_+\Bigl(\frac{1}{S(\overline\varphi)-S(\varphi_0)}\Bigr)^\frac{a(a+1)}{\theta} \\
            &\qquad + \tilde{n}^{-\frac{2}{\gamma+2}}\cdot(\log \tilde{n})^{\frac{2(2-\gamma)}{\theta(\gamma+2)}+\frac{2\gamma^\prime}{\gamma+2}}
            + \tilde{N}^{-\frac{2}{\gamma+2}}\cdot(\log \tilde{N})^{\frac{(2-\gamma)(a+2)(a+1)}{\theta(2+\gamma)}+\frac{2\gamma^\prime}{2+\gamma}}.
        \end{aligned}
    \end{equation}

    \textbf{(ii).} When $P$ is a polynomial-tailed distribution with $\EE_P\|X\|_2^{2a+4+c}<\infty$ for some $c>2a$,
    \begin{equation}
        \EE\|\nabla\hat\varphi_{n,N}-\nabla\varphi_0\|_{L^2(P)}^2\lesssim \Bigl(S(\overline\varphi)-S(\varphi_0)\Bigr)^{k} +  \tilde{n}^{-k_1}(\log \tilde{n})^{\gamma^\prime\cdot k_1}+ \tilde{N}^{-k_2}(\log \tilde{N})^{\gamma^\prime\cdot k_2},
    \end{equation}
    where $k=\frac{2+c}{2+c+(a+1)a}$, $k_1=\frac{2(2+c)}{(\gamma+2)c+4(a^2+a+2)}$ and $k_2=\frac{2(2+c)}{(\gamma+2)c+4a^2+(8-2\gamma)a+8}$.
\end{theorem}

\begin{remark}[Comments to Theorem \ref{thm: donsker function class with strongly convex, main text}]
    Ignoring the logarithmic factors, Theorem \ref{thm: donsker function class with strongly convex, main text} suggests that the estimation error for a sub-Weibull distribution $P$ is at the order of $\tilde{n}^{-2/(\gamma+2)}+\tilde{N}^{-2/(\gamma+2)}$ when $\gamma\leq 2$.
    This convergence rate is consistent with the findings in \cite{divol2022optimal}. Notably, for sub-Weibull distributions, the asymptotic convergence rate of the estimation error can potentially reach the parametric rate $\tilde{n}^{-1}+\tilde{N}^{-1}$ up to logarithmic terms, provided that $\Fcal$ is a Donsker class with $\gamma=0$, improving upon the results in Theorem \ref{thm: general function class with strongly convex, main text}.
    
    Additionally, for polynomial-tailed distributions, as $m$ increases, the convergence rate will approach $\tilde{n}^{-2/(\gamma+2)}+\tilde{N}^{-2/(\gamma+2)}$. The requirement of $c>2a$, rather than $c>0$, is due to the involvement of the discrete empirical measure $\QQ_N$. Specifically, for the semi-discrete estimator $\nabla\hat\varphi_n$ when $Q$ is known, assuming $\EE_P\|X\|_2^{2a+4+c}<\infty$ for some $c>0$ still allows for a convergence rate of $(S(\overline\varphi)-S(\varphi_0))^{k} +  \tilde{n}^{-k_1}(\log \tilde{n})^{\gamma^\prime\cdot k_1}$.
\end{remark}

\subsection{\texorpdfstring{$L^\infty$}{L-infty} Rates with Bounded Supports}\label{sec: L_infty Rates with Bounded Supports}

While we focus on $L^2(P)$ convergence rates, there has been a growing interest in studying $L^\infty$ convergence rates, as highlighted in \cite{ghosal2022multivariate}. In this section, we establish a link between the semi-local uniform convergence rates and the $L^2(P)$ rates through Proposition \ref{prop: semi-local uniform convergence rate, main text}. This result relies on the analytical properties of $(\beta,b)$-smooth functions, treating $\varphi_1$ and $\varphi_2$ as deterministic functions. 

\begin{proposition}\label{prop: semi-local uniform convergence rate, main text}
    Suppose probability measure $P$ has a density $p$ on $\RR^d$, such that the support $\Omega=\{x\in\RR^d:p(x)>0\}$ is closed, $\Omega^\circ$ is a Lipschitz domain (if $\Omega\not=\RR^d$) and the density $p$ is bounded from above. Assume $\varphi_1,\varphi_2\in C^2(\RR^d)$ are $(\beta,b)$-smooth, and $\nabla\varphi_1,\nabla\varphi_2$ are $L^2(P)$ integrable functions. Then for any $R>0$, we have 
    $$\sup_{x\in B(0,R)\cap\Omega}\|\nabla\varphi_1(x)-\nabla\varphi_2(x)\|_2 \lesssim \|\nabla\varphi_1-\nabla\varphi_2\|_{L^2(P)}^\frac{2}{d+2}\lor\|\nabla\varphi_1-\nabla\varphi_2\|_{L^2(P)},$$
    where the suppressed constant depends on $R$, $d$, $P$, $\beta$ and $b$.
\end{proposition}

Compared with Proposition F.1 in \cite{ghosal2022multivariate}, our result achieves the same rate with fewer and weaker assumptions. For instance, the support of probability measure $P$ can be unbounded, and $\nabla\varphi_0$ does not need to be locally uniformly Lipschitz. Though weaker than the uniform convergence rate, Proposition \ref{prop: semi-local uniform convergence rate, main text} is particularly helpful when $P$ has a bounded support. Combining this proposition with Theorems \ref{thm: general function class with strongly convex, main text} and \ref{thm: donsker function class with strongly convex, main text}, we obtain the following results.

\begin{corollary}\label{corollary: L_infty convergence rate for bounded support}
    Suppose Assumptions \ref{assumption: existence of Brenier's potential, main}, \ref{assumption: Convexity and Smoothness of Brenier's potential, main} and \ref{assumption: function class Fcal, main} hold. Additionally, assume the support $\Omega=\{x\in\RR^d:p(x)>0\}$ of $P$ is compact, $\Omega^\circ$ is a Lipschitz domain, and the density $p$ is bounded from above and below on $\Omega$. Then, for arbitrary $(\alpha,a)$-convex $\overline\varphi\in\Fcal$,
    
    \textbf{(i).} when $\Fcal$ is a Donsker function class: 
    $$\begin{aligned}
        \EE\|\nabla\hat\varphi_{n,N}-\nabla\varphi_0\|_{L^\infty(\Omega)} 
        \lesssim_{\log} \Bigl(S(\overline\varphi)-S(\varphi_0)\Bigr)^\frac{1}{2}\lor \Bigl(S(\overline\varphi)-S(\varphi_0)\Bigr)^\frac{1}{d+2}+ (\tilde{n}\land\tilde{N})^{-\frac{2}{(\gamma+2)(d+2)}}.
    \end{aligned}$$

    \textbf{(ii).} When $\Fcal$ is a general function class: 
    $$\begin{aligned}
        \EE\|\nabla\hat\varphi_{n,N}-\nabla\varphi_0\|_{L^\infty(\Omega)} 
        \lesssim_{\log} \Bigl(S(\overline\varphi)-S(\varphi_0)\Bigr)^\frac{1}{2}\lor \Bigl(S(\overline\varphi)-S(\varphi_0)\Bigr)^\frac{1}{d+2}+ (\tilde{n}\land\tilde{N})^{-\frac{1}{\gamma(d+2)}}.
    \end{aligned}$$
\end{corollary}

In Section \ref{sec: Our Poincaré-type Inequalities}, we will show that the assumptions in Corollary \ref{corollary: L_infty convergence rate for bounded support} serve as sufficient conditions for Assumption \ref{assumption: function class Fcal, main}.  Additionally, Proposition \ref{prop: semi-local uniform convergence rate, main text} can be applied to results in Sections \ref{sec: Relaxing the Strong Convexity} and \ref{sec: Beyond 4th-order Moments} as well, though these extensions are omitted here for simplicity.

\subsection{Relaxing the Strong Convexity}\label{sec: Relaxing the Strong Convexity}

\begin{table}[ht]
    \centering
    \begin{tabular}{c|c|c}
        \hline
        & \multicolumn{2}{c}{Brenier's potential $\varphi_0$ / OT map $\nabla\varphi_0$} \\ \hline
        & Convexity & Smoothness \\ \hline
        \cite{hutter2021minimax} & Strongly convex & $\alpha$-Hölder \\
        \cite{manole2021plugin} & Strongly convex & Lipschitz \\
        Theorem 2.2 in \cite{deb2021rates} & Convex & Lipschitz \\
        Theorem 2.5 in \cite{deb2021rates} & Convex & Lipschitz \\
        \cite{muzellec2021near} & Convex & $(m+2)$-Sobolev ($m>d+1$) \\
        \cite{gunsilius2022convergence} & Strictly Convex & Lipschitz \\
        \cite{divol2022optimal} & $(\alpha,a)$-convex  & $(\beta,a)$-smooth \\
        Our Theorems \ref{thm: general function class with strongly convex, main text}, \ref{thm: donsker function class with strongly convex, main text} & $(\alpha,a)$-convex  & $(\beta,a)$-smooth \\
        Our Theorems \ref{thm: general function class without strongly convex, main text}, \ref{thm: donsker function class without strongly convex, main text}, \ref{thm: convergence rate of semi-discrete optimal transport, second moment, main text} & Convex & $(\beta,b)$-smooth  \\ \hline
    \end{tabular}
    \caption{Assumptions on the Brenier's potential $\varphi_0$ / OT map $\nabla\varphi_0$ in the literature. ``Strong convexity'' of the Brenier's potential implies that its Hessian, $\nabla^2\varphi_0\succeq \alpha I_d$, is positive definite for some $\alpha>0$, and ``Lipschitz'' OT map means $\nabla\varphi_0$ is a Lipschitz function. Therefore, if the Brenier's potential is strongly convex and the corresponding OT map is Lipschitz, the OT map is bi-Lipschitz. The terms $(\beta,a)$-smooth and $(\beta,b)$-smooth refer to the Brenier's potential. In \cite{gunsilius2022convergence}, the Lipschitz property is implicitly assumed as per Theorem 10.4 in \cite{rockafellar1997convex}. }
    \label{table: assumptions of OT estimation, potential}
\end{table}

The combination of $(\alpha,a)$-convexity and $(\beta,a)$-smoothness assumptions of the Brenier's potential with the same $a\geq 0$ lays the foundation for establishing the convergence rates of the OT map estimator from Equation \eqref{eq: discrete-discrete OT} for both \cite{divol2022optimal} and our results in Theorems \ref{thm: general function class with strongly convex, main text} and \ref{thm: donsker function class with strongly convex, main text}. Indeed, this pair of assumptions generalizes those in existing literature, corresponding exactly to the bi-Lipschitz continuity of the OT map when $a=0$.

However, numerous situations arise where the Brenier's potential is not $(\alpha,a)$-convex for any $a\geq 0$. For instance, consider $P$ as the standard normal distribution $N(0,1)$ and $Q$ as the uniform measure on interval $(0,1)$. In this case, the Brenier's potential from $P$ to $Q$ is not $(\alpha,a)$-convex for any $a\geq 0$, even though the corresponding OT map (i.e. the CDF of $N(0,1)$) is Lipschitz. This discrepancy arises because the Hessian of the Brenier's potential, which corresponds to the normal density $\phi(x)=\exp(-x^2/2)/\sqrt{2\pi}$, is uniformly bounded from above by $+\infty$, but not bounded from below by $0$.

Beyond this intuitive example, Proposition \ref{prop: thickness and OT map convexity, smoothness} below, as a generalization of the reversed Caffarelli’s contraction theorem \citep{caffarelli2000monotonicity}, formalizes the relationship between the tail behavior of probability measures and the properties of the Brenier's potentials. Specifically, if the Brenier's potential $\varphi_0$ is $(\beta,b)$-smooth, then the tail of $Q$ cannot be much heavier than that of $P$. Conversely, if $\varphi_0$ is $(\alpha,a)$-convex, the tail of $P$ cannot be much lighter than that of $Q$.

When we assume the OT potentials to be both $(\alpha,a)$-convex and $(\beta,a)$-smooth for a common $a\geq 0$, this pair of assumptions imposes a narrow constraint on the relationship between the probability measures $P$ and $Q$ in terms of the thickness of their tails. Loosely speaking, under such assumptions, Proposition \ref{prop: thickness and OT map convexity, smoothness} implies that both $P$ and $Q$ should either be sub-Weibull measures, have polynomial tails, or process bounded supports. 

\begin{proposition}[Reversed Caffarelli’s contraction theorem]\label{prop: thickness and OT map convexity, smoothness}
    Let $\varphi_0$ be the Brenier's potential from $P$ to $Q$ on $\RR^d$.
    \begin{itemize}
        \item[\textit{\textbf{(i).}}] If $\varphi_0$ is a $(\beta,b)$-smooth for some $\beta>0$ and $b\geq 0$:
        \begin{itemize}
            \item[$\bullet$] $\EE_P[\exp(\|X\|_2^{\theta})]<\infty$ for some $\theta>0$ implies $\EE_Q[\exp(\|Y\|_2^{\theta/(b+1)})]<\infty$.
            \item[$\bullet$] $\EE_P\|X\|_2^{m}<\infty$ for some $m>0$ implies $\EE_Q\|Y\|_2^{m/(b+1)}<\infty$.  
        \end{itemize}

        \item[\textit{\textbf{(ii).}}] Conversely, if $\varphi_0$ is an $(\alpha,a)$-convex for some $\alpha>0$ and $a\geq 0$:
        \begin{itemize}
            \item[$\bullet$] $\EE_Q[\exp(\|Y\|_2^{\theta/(a+1)})]<\infty$ for some $\theta>0$ implies $\EE_P[\exp(\|X\|_2^{\theta})]<\infty$.
            \item[$\bullet$] $\EE_Q\|Y\|_2^{m/(a+1)}<\infty$ for some $m>0$ implies $\EE_P\|X\|_2^{m}<\infty$. 
        \end{itemize}
    \end{itemize}
\end{proposition}

To overcome the limitations introduced by the combination of $(\alpha,a)$-convexity and $(\beta,a)$-smoothness assumptions, a novel, sieve plug-in estimator of the Brenier's potential is introduced in Equation \eqref{eq: new discrete-discrete estimator, overview}. This new estimator allows us to establish the convergence rates without the $(\alpha,a)$-convexity assumption. Despite this advancement, our conclusions still hinge on the $(\beta,b)$-smoothness assumption. Therefore, it is natural to inquire about the limitations of making such an assumption. After all, the OT map from $U(0,1)$ to $N(0,1)$, i.e. the quantile function of $N(0,1)$, is not $(\beta,b)$-smooth for any $b\geq 0$. 

To address this issue, let $\varphi_{P\to Q}$ and $\varphi_{Q\to P}$ denote the Brenier's potentials from $P$ to $Q$, and from $Q$ to $P$, respectively. If $\varphi_{P\to Q}$ is not $(\beta,b)$-smooth for any $b\geq 0$, but $\varphi_{Q\to P}$ meets such assumption for some $b\geq 0$, we can first estimate the Brenier's potential from $Q$ to $P$ as $\hat\varphi_{Q\to P}$. Then we take the gradient of its convex conjugate as the estimated OT map from $P$ to $Q$: $\nabla\hat\varphi_{Q\to P}^\ast$. The validity of this approach is guaranteed by Proposition \ref{prop: regularity of OT map} below.

\begin{proposition}[Caffarelli’s regularity theory, Theorem 2.12 (iv) in \cite{villani2003topics}]\label{prop: regularity of OT map}
    Let $P$ and $Q$ be two probability measures on $\RR^d$ with finite second-order moments. Denote $\nabla\varphi_0$ as the OT map from $P$ to $Q$. If both $P$ and $Q$ are absolutely continuous with respect to the Lebesgue measure $\lambda_d$, then for $\dd P$ almost all $x$ and $\dd Q$ almost all $y$, $\nabla\varphi_0^\ast\circ\nabla\varphi_0(x)=x$, $\nabla\varphi_0\circ\nabla\varphi_0^\ast(y)=y$, and $\nabla\varphi_0^\ast$ is the OT map from $Q$ to $P$.
\end{proposition}

When the estimated OT map $\nabla\hat\varphi_{Q\to P}$ converges to $\nabla\varphi_{Q\to P}$ in $L^2(Q)$, the convergence of potential, $\hat\varphi_{Q\to P}^\ast\to \varphi_{P\to Q}$ in $L^2(P)$, can be ensured by Proposition \ref{prop: regularity of OT map}, Lemma 4 of \cite{gunsilius2022convergence} and the Poincaré-type inequalities in Section \ref{sec: Poincaré-type Inequalities}. Under some extra mild regularity conditions, we can apply Theorem 25.7 in \cite{rockafellar1997convex}, which asserts that the convexity of $\hat\varphi_{Q\to P}^\ast$ and $\varphi_{P\to Q}$ implies the convergence of $\nabla\hat\varphi_{Q\to P}^\ast$ to $\nabla\varphi_{P\to Q}$. For this reason, it suffices for us to focus only on Brenier's potentials that are $(\beta,b)$-smooth in this section:

\begin{assumption}[Smoothness of Brenier's potential]\label{assumption: Smoothness of Brenier's potential, main}
    Assume $\varphi_0$ is a $(\beta,b)$-smooth Brenier's potential from $P$ to $Q$ for some $\beta>0$ and $b\geq 0$. Function class $\Fcal$ is also made up of $(\beta,b)$-smooth functions with $\varphi(0)=0$ for all $\varphi\in\Fcal$. 
\end{assumption}

As the combination of $(\alpha,a)$-convexity and $(\beta,a)$-smoothness assumptions aims to control the empirical process involving $\QQ_N$ and the convex conjugates defined in Equation \eqref{eq, semi-dual of Kantorovich definition}, inspired by \cite{shen1994convergence}, we propose a ``sieve convex conjugate'' to replace it. Recall that for a positive and increasing sequence $\{M_n\}_n$, our sieve estimator is defined as 
\begin{equation}\label{eq: new discrete-discrete estimator}
    \tilde\varphi_{n,N}\in\argmin_{\varphi\in\Fcal}\PP_n\varphi+\QQ_N\varphi^{\ast,(n)},\quad\text{where}\quad \varphi^{\ast,(n)}(y)=\sup_{x\in B(0,M_n)}\inner{x,y}-\varphi(x).
\end{equation}

In this sieve convex conjugate, we restrict the supremum's search domain to a bounded hyperball $B(0,M_n)$, effectively serving as a pseudo-support for $P$. In fact, the first-order optimality condition tells that $\varphi^\ast(y)=\inner{(\nabla\varphi)^{-1}(y), y}-\varphi((\nabla\varphi)^{-1}(y))=\inner{\nabla\varphi^\ast(y), y}-\varphi(\nabla\varphi^\ast(y))$. By posing a restriction on the supremum's search domain, we are essentially restricting the range of $\nabla\varphi^\ast$, acting as an implicit penalty on the model complexity of convex conjugates. Subsequently, this approach allows us to bypass the necessity of the $(\alpha,a)$-convexity assumption and directly control the empirical process involving $\QQ_N$ and the sieve convex conjugates $\varphi^{\ast,(n)}$'s as follows:
\begin{lemma}\label{lemma: concentrtaion of QQ_N-Q for non-strongly convex EP, main text}
    Suppose Assumption \ref{assumption: function class Fcal, main} holds with $\gamma<2$, and $\Fcal$ is made up of $(\beta,b)$-smooth potentials. Then for any potential $\overline\varphi\in\Fcal$ and any function $h_n$ such that $0\leq h_n\leq 1$, $$\EE^\ast\Bigl[\sup_{\varphi\in\Fcal}|(Q-\QQ_N)(\varphi^{\ast,(n)}-\overline\varphi^{\ast,(n)})\cdot h_n|\Bigr]\lesssim \tilde{N}^{-\frac{1}{2}}\cdot M_n^{\frac{2-\gamma}{2}}+M_n\cdot\int_{\|y\|_2\geq N^{\frac{1}{2}}\cdot M_n^{\frac{\gamma}{2}}} \inner{y}\, Q(\dd y).$$
\end{lemma}
This lemma indicates that by appropriately selecting the sequence $\{M_n\}_n$ based on the tail thickness of $P$ and $Q$, we can control this empirical process without assuming the $(\alpha,a)$-convexity of $\varphi_0$. Additionally, a corresponding version of this lemma for a general function class (i.e. $\gamma\geq 2$) is available as well. 

For both theoretical analysis and practical implementations, we set $M_n=\max_{1\leq i\leq n}\|X_i\|_2$, representing the maximum norm among all samples in $\PP_n$. We hypothesize that alternative choices of $\{M_n\}_n$ or different shapes for the supremum's search domain could potentially yield better convergence rates, and we leave this exploration to interested readers.

We now present the results for general function class $\Fcal$ as a counterpart to Theorem \ref{thm: general function class with strongly convex, main text}:

\begin{theorem}\label{thm: general function class without strongly convex, main text}
    Suppose Assumptions \ref{assumption: existence of Brenier's potential, main}, \ref{assumption: function class Fcal, main} and \ref{assumption: Smoothness of Brenier's potential, main} hold. $\Fcal$ is a general function class in the sense that $\gamma\geq 2$. Let $\overline\varphi$ be an arbitrary potential in $\Fcal$.
    
    \textbf{(i).} Suppose $P$ is a sub-Weibull probability measure with parameters $(\theta,K)$ for some $\theta,K>0$. Then for large $n$ and $N$ such that $N\lesssim n$,
    \begin{equation}
        \begin{aligned}
            &\EE\|\nabla\tilde\varphi_{n,N}-\nabla\varphi_0\|_{L^2(P)}^2
            \lesssim_{\log\log}\Bigl(S(\overline\varphi)-S(\varphi_0)\Bigr)\cdot\log_+\Bigl(\frac{1}{S(\overline\varphi)-S(\varphi_0)}\Bigr)^\frac{b(b+1)}{\theta}\\
            &\qquad+\tilde{n}^{-\frac{1}{\gamma}}\cdot(\log\tilde{n})^{b+2+\frac{\gamma^\prime\lor 1}{2}+\frac{b(b+1)}{\theta}}
            +\tilde{N}^{-\frac{1}{\gamma}}\cdot(\log \tilde{N})^{1+\frac{\gamma^\prime\lor 1}{2}+\frac{b(b+1)}{\theta}}\cdot (\log n)^{1/\theta}.
        \end{aligned}
    \end{equation}

    \textbf{(ii).} Suppose $P$ is a polynomial-tailed distribution with $\EE_P\|X\|_2^m<\infty$ for some $m>2b+4$. Then for large $n$ and $N$ such that $n^{\frac{1}{m}\bigl(\frac{1}{\gamma}\land \frac{m-(2b+2)}{2m}\bigr)}\lesssim N\lesssim n$,
    \begin{equation}
        \begin{aligned}
            &\EE\|\nabla\tilde\varphi_{n,N}-\nabla\varphi_0\|_{L^2(P)}^2
            \lesssim\Bigl(S(\overline\varphi)-S(\varphi_0)\Bigr)^{\frac{m-2(b+1)}{m+(b+1)(b-2)}}\\
            &\qquad+\tilde{n}^{-(\frac{1}{\gamma}\land\frac{m-(2b+4)}{2m})\frac{m-2(b+1)}{m+(b+1)(b-2)}}\cdot(\log \tilde{n})^{\frac{\gamma^\prime+2}{2}\cdot\frac{m-2(b+1)}{m+(b+1)(b-2)}}+\tilde{n}^{\frac{b+2-m}{m}\frac{m-2(b+1)}{m+(b+1)(b-2)}}\\
            &\qquad+\tilde{N}^{-(\frac{1}{\gamma}\land\frac{m-(2b+2)}{2m})\frac{m-2(b+1)}{m+(b+1)(b-2)}}\cdot n^{\frac{m-2(b+1)}{m^2+m(b+1)(b-2)}}\cdot (\log N)^{\frac{2+\gamma^\prime\lor 1}{2}\frac{m-2(b+1)}{m+(b+1)(b-2)}}.\\
        \end{aligned}
    \end{equation}
\end{theorem}

Similarly, when $\Fcal$ is a Donsker class, Theorem \ref{thm: donsker function class with strongly convex, main text} can also be generalized as follows:

\begin{theorem}\label{thm: donsker function class without strongly convex, main text}
    Suppose Assumptions \ref{assumption: existence of Brenier's potential, main}, \ref{assumption: function class Fcal, main}, \ref{assumption: Poincare-type inequalities, main text} and \ref{assumption: Smoothness of Brenier's potential, main} hold. $\Fcal$ is a Donsker function class in the sense that $\gamma<2$. Let $\overline\varphi$ be an arbitrary potential in $\Fcal$.

    \textbf{(i).} When $P$ is a sub-Weibull probability measure with parameters $(\theta,K)$ for some $\theta,K>0$. Then, for large $n$ and $N$ such that $N\lesssim n$,
    \begin{equation}
        \begin{aligned}
            &\EE\|\nabla\tilde\varphi_{n,N}-\nabla\varphi_0\|_{L^2(P)}^2
            \lesssim_{\log\log} \Bigl(S(\overline\varphi)-S(\varphi_0)\Bigr)\cdot\log_+\Bigl(\frac{1}{S(\overline\varphi)-S(\varphi_0)}\Bigr)^\frac{b(b+1)}{\theta} \\
            &\qquad+ \tilde{n}^{-\frac{2}{\gamma+2}}\cdot(\log \tilde{n})^{\frac{2(2-\gamma)}{\theta(\gamma+2)}+\frac{2\gamma^\prime}{\gamma+2}}
            + \tilde{N}^{-\frac{1}{2}}\cdot (\log n)^{\frac{2-\gamma}{2\theta}}\cdot (\log \tilde{N})^{\frac{b(b+1)}{\theta}}.
        \end{aligned}
    \end{equation}

    \textbf{(ii).} When $P$ is a polynomial-tailed distribution with $\EE_P\|X\|_2^{2b+4+c}<\infty$ for some $c>0$. Then, for large $n$ and $N$ such that $n^\frac{2-\gamma}{2b+4+c}\lesssim N\lesssim n$,
    \begin{equation}
        \EE\|\nabla\tilde\varphi_{n,N}-\nabla\varphi_0\|_{L^2(P)}^2\lesssim \Bigl(S(\overline\varphi)-S(\varphi_0)\Bigr)^{k} + \tilde{n}^{-k_1}(\log \tilde{n})^{\gamma^\prime\cdot k_1}+ \tilde{N}^{-\frac{k}{2}}\cdot \tilde{n}^{\frac{(2-\gamma)k}{2(2b+4+c)}}+\tilde{n}^{-k_2},
    \end{equation}
    where $$k=\frac{2+c}{2+c+(b+1)b};\qquad k_1=\frac{2(2+c)}{(\gamma+2)c+4(b^2+b+2)};\qquad k_2=k\cdot \frac{b+2+c}{2b+4+c}.$$
\end{theorem}

\begin{remark}[Comments to Theorem \ref{thm: donsker function class without strongly convex, main text}]
    Different from Theorem \ref{thm: donsker function class with strongly convex, main text} that the estimation error is at the order of $\tilde{n}^{-\frac{2}{\gamma+2}}+\tilde{N}^{-\frac{2}{\gamma+2}}$ for sub-Weibull $P$, Theorem \ref{thm: donsker function class without strongly convex, main text} established a slower estimation error of $\tilde{n}^{-\frac{2}{\gamma+2}}+\tilde{N}^{-\frac{1}{2}}$. This difference arises because Theorem \ref{thm: donsker function class with strongly convex, main text} is established with the uniform localization technique in Theorem 3.4.1 of \cite{vaart2023empirical}. However, we are unable to provide uniform localized convergence rate in Lemma \ref{lemma: concentrtaion of QQ_N-Q for non-strongly convex EP, main text}. Despite the slower theoretical convergence rate, the second estimator demonstrates superior performance and robustness in numerical experiments, as evidenced in Section \ref{sec: Simulation Results}.
\end{remark}

\subsection{Relaxing Fourth-order Moments}\label{sec: Beyond 4th-order Moments}

Although we have derived the convergence rates for sub-Weibull and polynomial-tailed distributions and removed the $(\alpha,a)$-convexity assumption on the Brenier's potential, we still assume $P$ processes fourth-order moments. This requirement arises because the existing empirical process tools we use (Theorem 16 in \cite{von2004distance} and Remark 3.5.14 in \cite{gine2021mathematical}) necessitate that the function class $\Fcal$ consists of square-integrable functions. However, to ensure the quantities $\|\nabla\varphi-\nabla\varphi_0\|_{L^2(P)}^2$ and $P\varphi+Q\varphi^\ast$ exist for $(\beta,b)$-smooth potentials, it suffices to assume that $\EE_P\|X\|_2^{2b+2}<\infty$. Moreover, our assumption that $P$ has fourth-order moments is stricter than the assumption in Brenier's Theorem, which only requires $P$ and $Q$ to have second-order moments.

Thus, it is natural to ask whether our convergence results can be further extended to approach such theoretical limit in the context of Brenier's Theorem? To address this question, given that $(\beta,b)$-smooth potentials are not $L^2(P)$ integrable under weaker moment conditions, we have developed the following maximal inequality:

\begin{lemma}\label{lemma: convergence of EP with L^1 integrable functions}
    Given a function class $\Fcal$, denote its envelope function $F(x)=B\inner{x}^{b+2}$ for some $B>0$, $b\geq 0$. Assume $\EE_P\|X\|_2^{b+2+c}<\infty$ for some $c>0$, and there exists some $D_\Fcal>1$ and a function $H(x):(0,\infty)\to[0,\infty)$ such that any $x>0$, $$\sup_{P^\prime\in P_n}\log\Ncal(x,\Fcal,L^1(P^\prime))\leq D_\Fcal\cdot H(x)\quad\text{ and }\quad\int_0^{\infty}\sqrt{\frac{H(x)}{x}}\, \dd x<\infty\quad\text{ is integrable,}$$ where the supremum is taken over the product measure $P_n=P^{\otimes n}$. Then, $$\EE^\ast\|\PP_n-P\|_\Fcal\lesssim D_\Fcal^\frac{c}{2(b+2+c)}\cdot n^{-\frac{c^2}{2(b+2+c)^2}},$$ where the suppressed constant depends on $B$, $b$, $c$ and $\EE_P\|X\|_2^{b+2+c}$.
\end{lemma}

Now based on Lemma \ref{lemma: convergence of EP with L^1 integrable functions},
we are ready to present the convergence rates of the semi-discrete estimator $\nabla\hat\varphi_n$ and the new discrete-discrete estimator $\nabla\tilde\varphi_{n,N}$ in Theorem \ref{thm: convergence rate of semi-discrete optimal transport, second moment, main text}. 

\begin{theorem}\label{thm: convergence rate of semi-discrete optimal transport, second moment, main text}
    Suppose Assumptions \ref{assumption: existence of Brenier's potential, main}, \ref{assumption: Smoothness of Brenier's potential, main} and Lemma \ref{lemma: convergence of EP with L^1 integrable functions} hold. $\EE_P\|X\|_2^{b+2+c}<\infty$ for some $c>b$. $\Fcal$ is a Donsker class with $\gamma<1$. Let $\overline\varphi$ be an arbitrary potential in $\Fcal$.
    
    \textbf{(i).} For large $n$, the semi-discrete OT map estimator $\hat\varphi_{n}$ in Equation \eqref{eq: semi-discrete OT} satisfies
    \begin{equation}
        \EE\|\nabla\hat\varphi_{n}-\nabla\varphi_0\|_{L^2(P)}^2\lesssim \Bigl(S(\overline\varphi)-S(\varphi_0)\Bigr)^{\frac{c-b}{c+b^2}} + D_\Fcal^\frac{c(c-b)}{2(b+2+c)(c+b^2)}\cdot n^{-\frac{c^2(c-b)}{2(b+2+c)^2(c+b^2)}}.
    \end{equation}

    \textbf{(ii).} For large $n$ and $N$ such that $n^\frac{2-\gamma}{b+2+c}\lesssim N\lesssim n$, the new discrete-discrete estimator in Equation \eqref{eq: new discrete-discrete estimator, overview} satisfies
    \begin{equation}
        \begin{aligned}
            \EE\|\nabla\tilde\varphi_{n,N}-\nabla\varphi_0\|_{L^2(P)}^2
            \lesssim& \Bigl(S(\overline\varphi)-S(\varphi_0)\Bigr)^{\frac{c-b}{c+b^2}}
             + D_\Fcal^\frac{c(c-b)}{2(b+2+c)(c+b^2)}\cdot n^{-\frac{c^2(c-b)}{2(b+2+c)^2(c+b^2)}}\\
            &+\tilde{N}^{-\frac{c-b}{2(c+b^2)}}\cdot \tilde{n}^{\frac{(2-\gamma)(c-b)}{2(b+2+c)(c+b^2)}}+\tilde{n}^{-\frac{c(c-b)}{(b+2+c)(c+b^2)}}.
        \end{aligned}
    \end{equation}
\end{theorem}

\begin{remark}[Comments on the assumption $c>b$ in Theorem \ref{thm: convergence rate of semi-discrete optimal transport, second moment, main text}]
    Despite the moment condition $\EE_P\|X\|_2^{b+2}<\infty$ suffices to ensure that $(\beta,b)$-smooth potentials $\varphi$'s are $L^1(P)$-integrable, we assume $\EE_P\|X\|_2^{b+2+c}<\infty$ for some $c>b$ for technical reasons in the proof. A rationale behind this is that, according to Proposition \ref{prop: thickness and OT map convexity, smoothness}, the probability measure $Q$ must satisfy $\EE_Q\|Y\|_2^{\frac{b+2+c}{b+1}}<\infty$. Thus, the condition $c>b$ ensures that $\frac{b+2+c}{b+1}>2$, meaning that $Q$ has finite second-order moments --- an assumption required by Brenier's Theorem. 
\end{remark}

Note that, unlike our previous results, the convergence rates in Theorem \ref{thm: convergence rate of semi-discrete optimal transport, second moment, main text} are established only when $\Fcal$ satisfies $\gamma<1$. That is because the integrability of $\sqrt{H(x)/x}=D_\Fcal^{\frac{1}{2}}\cdot x^{-\frac{\gamma+1}{2}}\cdot \log_+(1/x)^{\gamma^\prime/2}$ around $0$ holds only if $\gamma<1$.

\section{Poincaré-type Inequalities}\label{sec: Poincaré-type Inequalities}
As a distinctive and pivotal assumption in \cite{gunsilius2022convergence} and \cite{divol2022optimal}, the (global) Poincaré inequality plays a pivotal role in the recent study of optimal transport. This inequality, fundamental in PDE and functional analysis, implies the existence of a constant $C>0$, dependent solely on measure $P$, such that for any Sobolev function $u\in W^{1,2}(P)=H^1(P)$, the variance of $u$ can be controlled by the $L^2(P)$ norm of its gradient:
\begin{equation}\label{eq: poincare inequality}
    \var_P(u)\leq C\cdot \|\nabla u\|_{L^2(P)}^2.
\end{equation}

In \cite{gunsilius2022convergence}, the significance of the Poincaré inequality emerges in establishing convergence rates for optimal transport estimators using empirical process theory. Specifically, the study demonstrated that the minimum value of Equation \eqref{eq: discrete-discrete OT} can be well-separated, a key requirement in Theorem 3.2.5 of \cite{vaart2023empirical}, with the help of the Poincaré inequality. Despite the introduction of the Poincaré inequality does not yield stronger conclusions compared to others, it underscores its role in regularizing the behavior of the objective function. Subsequently, \cite{divol2022optimal} leveraged the Poincaré inequality to extend prior results to unbounded probability measures.

While extending convergence rates of OT map estimators to unbounded probability measures is a significant improvement, it is important to note that findings in \cite{divol2022optimal} are not universally applicable to all probability measures. Specifically, if a probability measure satisfies the Poincaré inequality, it must be sub-exponential \citep{bobkov1997poincare}. However, Brenier's theorem holds for probability measures with second-order moments, illustrating a notable distinction between this theoretical limit and the sub-exponential constraint. Thus, to bridge this gap, alternatives to the original Poincaré inequality need to be explored.

In the remainder of this section, we begin by reviewing the Muckenhoupt condition from harmonic analysis in Section \ref{sec: Local Muckenhoupt Condition}, which serves as a sufficient condition to the Poincaré inequality. Based on this, we hypothesize that establishing local Poincaré inequalities for probability measures beyond sub-exponential might be feasible. In Section \ref{sec: Our Poincaré-type Inequalities}, we introduce additional regularity conditions to develop our tailored Poincaré-type inequalities, presented in Assumption \ref{assumption: Poincare-type inequalities, main text}. These conditions stem from both direct consequences and commonly accepted sufficient conditions related to the global Poincaré inequality.

\subsection{Local Muckenhoupt Condition}\label{sec: Local Muckenhoupt Condition}

Before presenting our Poincaré-type inequalities, it is essential to review the origin of the global Poincaré inequality, which serves as the motivation for our results. The global Poincaré inequality arises from the Hardy-Littlewood maximal inequality \citep{chanillo1985weighted}, which is equivalent to the Muckenhoupt condition \citep{muckenhoupt1972weighted}. Namely, a probability density $p(\cdot)$ satisfies the Muckenhoupt condition if there exists a constant $C>0$, such that for any ball $B=B(x_0,r)\subset\RR^d$,
\begin{equation}\label{eq: muckenhoupt condition}
    \Bigl(\frac{1}{\lambda_d(B)}\int_B p(x)\, \dd x\Bigr)\Bigl(\frac{1}{\lambda_d(B)}\int_B \frac{1}{p(x)}\, \dd x\Bigr)\leq C<\infty,
\end{equation}
where $\lambda_d(B)$ is the Lebesgue measure of the ball $B$ in $\RR^d$.

This condition ensures that the density $p(\cdot)$ is neither too concentrated nor too dispersed within any given ball, thereby controlling the oscillation of $p(\cdot)$ and providing a foundation for inequalities like the Poincaré inequality. 

Leveraging this insight, we explore the possibility of establishing local Poincaré inequalities for probability measures extending beyond the sub-exponential class by verifying the Muckenhoupt condition locally. For simplicity, consider a probability measure $P$ with positive density function $p(\cdot)$ defined on the entire $\RR^d$. Suppose the logarithm of the density function, $\log p(\cdot)$, is Lipschitz continuous, i.e., there exists a constant $M>0$ such that for any $x_1,x_2\in\RR^d$,
\begin{equation}\label{eq: Lipschitz log-density}
    |\log p(x_1)-\log p(x_2)|\leq M\cdot \|x_1-x_2\|_2.
\end{equation}  
The Lipschitz continuity of the logarithm of the density function implies that $P$ has a tail heavier than that of the exponential distribution. Specifically, by setting $x_1=x$ and $x_2 = 0$ in Equation \eqref{eq: Lipschitz log-density}, we obtain $p(x)\geq p(0) e^{-M\cdot\|x\|_2}$. 

Moreover, for any ball $B = B(x_0, r)$ centered at $x_0$ with radius $r$, and any $x\in B$, we can derive from Equation \eqref{eq: Lipschitz log-density} that $$e^{-Mr}\leq e^{-M\cdot \|x-x_0\|_2}\leq p(x)/p(x_0)\leq e^{M\cdot \|x-x_0\|_2}\leq e^{M r}.$$ Substituting these upper and lower bounds into the left-hand side of Equation \eqref{eq: muckenhoupt condition} yields:
\begin{equation}\label{eq: local muckenhoupt condition}
    \Bigl(\frac{1}{\lambda_d(B)}\int_B p(x)\, \dd x\Bigr)\Bigl(\frac{1}{\lambda_d(B)}\int_B \frac{1}{p(x)}\, \dd x\Bigr)\leq e^{2Mr}.
\end{equation}

We observe that although the Muckenhoupt condition in Equation \eqref{eq: muckenhoupt condition} is not met by all balls in $\RR^d$, it is satisfied by those with uniformly bounded radii. Therefore, we hypothesize that it is possible for a probability measure beyond sub-exponential class to satisfy the Muckenhoupt condition locally, leading to a local Poincaré inequality as a consequence.

\subsection{Our Poincaré-type Inequalities}\label{sec: Our Poincaré-type Inequalities}

Although the global Muckenhoupt condition in Equation \eqref{eq: muckenhoupt condition} supports the global Poincaré inequality in Equation \eqref{eq: poincare inequality}, Equation \eqref{eq: local muckenhoupt condition} arises in a local sense. Regrettably, to the best of our knowledge, we only know that the local Poincaré inequality can be established by local Muckenhoupt condition in one-dimensional $\RR$ \citep{bjorn2020locally}. Consequently, to derive the Poincaré-type inequalities for probability measures beyond sub-exponential class in our study, we must pursue an alternative approach by introducing additional regularity conditions on the probability measure $P$, inspired by common sufficient conditions and consequences related to the global Poincaré inequality.

The first requirement we consider is the connectedness of the support $\Omega$ of probability measure $P$. Notably, if $P$ supports the global Poincaré inequality, $\Omega$ must be connected (Proposition 8.1.6 in \cite{heinonen2015sobolev}). Moreover, the global Poincaré inequality implies $\Omega$ must be quasi-convex (Theorem 8.3.2 in \cite{heinonen2015sobolev}), meaning that any two points in $\Omega$ can be connected via a curve whose length is at most proportional to the Euclidean distance between these two points.
For simplicity, we consider the sufficient condition of quasi-convexity that $\Omega$ is closed (Theorem 2.12, Theorem 2.10 in \cite{heinonen2005lectures} and McShane’s Extension Theorem \citep{mcshane1934extension}).

Lastly, as the doubling measure is a common condition for verifying the Poincaré inequality, we also require our probability measure $P$ exhibits a doubling property. To achieve it, we impose that the boundary of the support $\Omega$ is Lipschitz. Additionally, the probability density should be locally bounded away from below by $0$ and from above by $+\infty$.

Given these regularity conditions on the topological properties of $\Omega$ and the local boundedness of the density, and focusing on the collection of $(\beta,b)$-smooth potentials, we now present our Poincaré-type inequalities for sub-Weibull and polynomial-tailed distributions. 

\begin{assumption}\label{assumption: measure P for Poincare-type inequality}
    Suppose $P$ is a probability measure on $\RR^d$ with a density $p(x)$ with respect to Lebesgue measure. Denote the support as $\Omega=\{x\in\RR^d:p(x)>0\}$. Assume that
    \begin{enumerate}
        \item[\textit{\textbf{(1).}}] $\Omega$ is closed, connected and locally connected;
        \item[\textit{\textbf{(2).}}] if $\Omega\not=\RR^d$, the interior $\Omega^\circ$ is a Lipschitz domain;
        \item[\textit{\textbf{(3).}}] $p(x)$ is locally bounded on $\Omega$; \footnote{We say that the density $p(x)$ is locally bounded on $\Omega$, if for any $x\in \Omega$, there exist constants $r>0$ and $0<m\leq M<\infty$ such that $m\leq p(y)\leq M$ for all $y\in\Omega$ with $\|x-y\|_2\leq r$.}
    \end{enumerate}
\end{assumption}

\begin{proposition}[Poincaré-type Inequalities]\label{proposition: poincare-type inequality}
    Under Assumption \ref{assumption: measure P for Poincare-type inequality}, suppose $\varphi_1$ and $\varphi_2$ are two $(\beta,b)$-smooth functions such that $\varphi_1(0)=\varphi_2(0)=0$. Then, 

    \textbf{(i).} If $P$ is sub-Weibull with parameter $\theta$, we have 
    \begin{equation}\label{eq: poincare-type inequality for sub-Weibull distributions}
        \var_P(\varphi_1-\varphi_2)
        \lesssim \|\nabla\varphi_1-\nabla\varphi_2\|_{L^2(P)}^2\cdot \log_+\Bigl(\frac{1}{\|\nabla\varphi_1-\nabla\varphi_2\|_{L^2(P)}}\Bigr)^{2/\theta}.
    \end{equation} 
    
    \textbf{(ii).} If $\EE_P\|X\|_2^{2b+4+c}<\infty$ for some $c>0$, we have 
    \begin{equation}\label{eq: poincare-type inequality for polynomial-tailed distributions}
        \var_P(\varphi_1-\varphi_2)\lesssim \|\nabla\varphi_1-\nabla\varphi_2\|_{L^2(P)}^2\lor \|\nabla\varphi_1-\nabla\varphi_2\|_{L^2(P)}^\frac{2c}{c+2}.
    \end{equation}

    \textbf{(iii).} If $\EE_P\|X\|_2^{b+2+c}<\infty$ for some $c>0$, we have 
    \begin{equation}\label{eq: poincare-type inequality for polynomial-tailed distributions with second moment}
        \EE_P\Bigl|(\varphi_1-\varphi_2)-\EE_P[(\varphi_1-\varphi_2)]\Bigr|
        \lesssim \|\nabla\varphi_1-\nabla\varphi_2\|_{L^2(P)}\lor \|\nabla\varphi_1-\nabla\varphi_2\|_{L^2(P)}^\frac{c}{c+1}.
    \end{equation}
    
    For all the three inequalities above, the suppressed constants depend only on measure $P$, coefficients $\beta$, $b$ and $c$, and do not depend on functions $\varphi_1$ and $\varphi_2$.
\end{proposition}

\begin{remark}[Comments to Proposition \ref{proposition: poincare-type inequality}]
    A noteworthy observation is that the exponent $\frac{2c}{c+2}$ in Equation \eqref{eq: poincare-type inequality for polynomial-tailed distributions} approaches $2$ as $c\to\infty$. This suggests that as higher-order moments of $P$ exist, Equation \eqref{eq: poincare-type inequality for polynomial-tailed distributions} increasingly resembles the global Poincaré inequality. Meanwhile, the discrepancy between Equation \eqref{eq: poincare-type inequality for sub-Weibull distributions} and the global Poincaré inequality is confined to an additional logarithmic term. Therefore, we should expect similar convergence rates of plug-in OT map estimators for sub-Weibull measures as in \cite{divol2022optimal}. 

    Although Lemma \ref{lemma: convergence of EP with L^1 integrable functions} is established independently of the Poincaré-type inequality, we still present the corresponding $(1,2)$ Poincaré-type inequality in Equation \eqref{eq: poincare-type inequality for polynomial-tailed distributions with second moment}, which could be particularly useful if Theorem 3.5.13 in \cite{gine2021mathematical} can be extended to the $L^1(P)$ function class.
\end{remark}

\begin{remark}[Comments to the new assumptions in \cite{divol2022optimal}]
    It is worth noting that, concurrent with our work, the second preprint version of \cite{divol2022optimal} on \texttt{ArXiv} introduced similar assumptions regarding the local Poincaré inequality (Assumption C3) and locally doubling measure (Assumption C4). However, in their framework, these assumptions primarily serve as sufficient conditions for their Assumption C2 (see their Lemma 1), which refines the covering entropy of the function class $\Fcal$ in Equation \eqref{eq: discrete-discrete OT}. They still rely on the global Poincaré inequality to establish the convergence of empirical processes. Moreover, it can be demonstrated that our mild assumptions in Proposition \ref{proposition: poincare-type inequality} are also sufficient conditions for their Assumptions C3 and C4, as shown in our supplementary material (see \cite{ding2024supplement}). 
\end{remark}

\section{Numerical Experiments}\label{sec: numerical experiments}

In Section \ref{sec: Practical Consideration}, we review existing OT map estimators from a computational perspective. Traditional methods often face scalability challenges related to both sample size and dimensionality. In contrast, the machine learning community has made significant strides by leveraging modern neural networks and the semi-dual form of the Kantorovich problem as objective functions. Building on these advancements, we present our algorithm at the end of Section \ref{sec: Practical Consideration}. Subsequently, in Section \ref{sec: Simulation Results}, we provide numerical results demonstrating that the strong representation capabilities of neural networks facilitate the estimation of OT maps in a time-efficient, resource-effective, and scalable manner.

\subsection{Practical Consideration}\label{sec: Practical Consideration}

\cite{hutter2021minimax} explored three different estimators in Section 6.1 of their study. The first baseline estimator (in their Section 6.1.1) is the conditional mean of coupling distribution given $X_i$'s, which is the solution to the Kantorovich problem (Equation \eqref{eq: kantorovich problem definition}) in its empirical form. As the coupling is the solution to a linear program, this estimator suffers from computational complexity, limiting its scalability with increasing sample size. Additionally, the empirical Kantorovich problem restricts the estimator from extrapolating beyond the observations in $\PP_n$. \cite{deb2021rates} presents an estimator resembling this baseline estimator, involving the replacement of the empirical measure with kernel density estimations of $P$ and $Q$ in the Kantorovich problem. However, this innovation restricts itself to probability measures with bounded domains and smooth densities.

The second estimator in Section 6.1.2 of \cite{hutter2021minimax} employs wavelets and requires discretizing the domain of measures $P$ and $Q$ into grid points, making it impractical for high-dimensional spaces. Building on this, \cite{manole2021plugin} proposed a minimax-optimal wavelet estimator under additional smoothness assumptions on the probability densities. However, this approach still faces significant computational challenges, highlighting a substantial statistical-computational gap.

The third estimator in Section 6.1.3 of \cite{hutter2021minimax} utilizes Reproducing Kernel Hilbert Spaces (RKHS) with Gaussian radial basis functions. Its objective function is formulated based on matching observed samples $X_i$'s and $Y_j$'s after transformation. This matching approach requires an equal number of observations from measures $P$ and $Q$ (i.e. $n=N$). Inspired by nonparametric least squares regression, \cite{manole2021plugin} extended this estimator to accommodate different sample sizes by proposing the "Convex Least Squares Estimator," which reduces to solving a finite-dimensional quadratic program. Nevertheless, the empirical Kantorovich problem remains a computational bottleneck, similar to the baseline estimator.

\cite{gunsilius2022convergence} addressed some of these challenges by utilizing the semi-dual form of the Kantorovich problem. Although their estimator exhibits a slower convergence rate compared to \cite{hutter2021minimax} due to the replacement of empirical measures with smooth kernel estimations in Equation \eqref{eq: discrete-discrete OT}, it offers simpler implementation and greater conceptual clarity than the theoretical wavelet estimator. Building on this approach, \cite{divol2022optimal} further analyzed the estimator using the original discrete empirical measures $\PP_n$ and $\QQ_n$, extending convergence rate studies to more general settings as discussed in Section \ref{sec: main results}. Inspired by the sieve method, our new estimator in Equation \eqref{eq: new discrete-discrete estimator, overview} further relaxes the technical $(\alpha,a)$-convexity assumption on Brenier's potential.

Although the OT map estimators derived from the semi-dual form of the Kantorovich problem possess excellent theoretical properties, the convex conjugate poses challenges for convenient computation. Thankfully, the semi-dual form of the Kantorovich problem has the potential to bridge the statistical-computational gap. Relying on the sum-of-squares (SoS) tight reformulation of OT proposed by \cite{vacher2021dimension}, \cite{muzellec2021near} converted the semi-dual form of the Kantorovich problem into an unconstrained convex program and proved the near-optimality of their estimator. However, their transformation is built on a strong smoothness assumption requiring both Brenier's potentials $\varphi_0$ and its convex conjugate $\varphi_0^\ast$ belong to the Sobolev space $H^{m+2}(\RR^d)$ with $m>d+1$ (see their Theorem 1). Moreover, their empirical optimization program (see their Equation (5)) is designed for equal numbers of observations from measures $P$ and $Q$. Despite the merit of bridging the statistical-computational gap, such limitations are undesirable in terms of removing as many constraints on Brenier's potential as possible.

Meanwhile, substantial progress in computing the OT map estimators has been made in the machine learning community. Leveraging the fact that $\inner{y,\nabla \psi(y)}-\varphi(\nabla\psi(y))\leq\varphi^\ast(y)$ for all function $\psi$, \cite{makkuva2020optimal} proposed the following minimax optimization program $$\min_{\varphi\in\Fcal}\max_{\psi\in\Fcal}\EE_P[\varphi(X)]+\EE_Q[\inner{Y,\nabla\psi(Y)}-\varphi(\nabla\psi(Y))].$$ To address the optimization difficulties inherent in the minimax optimization, \cite{korotin2019wasserstein} replaced the maximization with the cycle-consistency loss from CycleGAN \citep{zhu2017unpaired}. Incorporating a hyperparameter $\lambda>0$, the objective function becomes $$\min_{\varphi,\psi\in\Fcal}\EE_P[\varphi(X)]+\EE_Q[\inner{Y,\nabla\psi(Y)}-\varphi(\nabla\psi(Y))]+\lambda\cdot\EE_Q[\|\nabla\varphi\circ\nabla\psi(Y)-Y\|_2^2].$$ Ideally, $\psi$ should be the convex conjugate of the Brenier's potential, making the cycle-consistency loss a variant of Proposition \ref{prop: regularity of OT map}. 

Different from \cite{makkuva2020optimal} and \cite{korotin2019wasserstein} that utilize two convex potentials $\varphi$ and $\psi$, \cite{huang2020convex} adopted the semi-dual form of the Kantorovich problem in Equation \eqref{eq, semi-dual of Kantorovich definition} as the objective function directly. Specifically, given a convex potential $\varphi$, a numerical optimizer solves for its corresponding convex conjugate. Since maximizing $x\mapsto \inner{x,y}-\varphi(x)$ can be reformulated as a convex optimization problem, computing the convex conjugate of $\varphi$ is straightforward and computationally efficient. Furthermore, although \cite{huang2020convex} requires $\varphi$ to be strongly convex by adding a quadratic term: $\varphi_\alpha(x)=\varphi(x)+\frac{\alpha}{2}\|x\|_2^2$, our experiments indicate that this optimization method works without this quadratic term as well.

Beyond innovations in objective functions, these machine learning works (e.g., \cite{makkuva2020optimal,korotin2019wasserstein,huang2020convex}) employ input convex neural network (ICNN) from \cite{amos2017input} as their working function class. Unlike traditional nonparametric models, neural networks exhibit powerful capabilities in function estimation and complex data modeling. As a class of neural network specifically designed for fitting convex functions, ICNNs share some structural similarities with ordinary fully-connected neural networks. Denote $L(\cdot)$ as a linear layer and $L^+(\cdot)$ as a linear layer with non-negative weights, an ICNN with $K$ layers can be defined recursively as $$\varphi(x)=L_K^+(s(z_{K-1}))+L_K(x);\quad z_k:=L_k^+(s(z_{k-1}))+L_k(x);\qquad z_1=L_1(x),$$ where $s(\cdot)$ is a non-decreasing, convex activation function. By utilizing ICNNs, it is possible to achieve a scalable, efficient and end-to-end approach to OT map estimation.

\begin{algorithm}
    \caption{Plug-in OT Map Estimators $\nabla\hat\varphi_{n,N}$ and $\nabla\tilde\varphi_{n,N}$ with ICNN}
    \label{alg:ot_map_icnn}
    \begin{algorithmic}[1]
        \Require Samples $(X_i)_{i=1}^n$ and $(Y_j)_{j=1}^N$; ICNN model $\Fcal_\Theta$; number of epochs $T$; batch sizes $m, M$.
        \State Initialize $\varphi_\theta\in\Fcal_\Theta$
        \For{$t=1, \cdots, T$}
            \For{mini-batch $(X_{i_k})_{k=1}^m$ in $(X_i)_{i=1}^n$, and $(Y_{j_k})_{j=1}^M$ in $(Y_j)_{j=1}^N$}
                \State Compute $(\varphi_\theta^\ast(Y_{j_k}))_{k=1}^M$ with Algorithm~\ref{alg:convex_conjugate}
                \State Compute loss: $L\leftarrow\frac{1}{m}\sum_{k=1}^m \varphi_\theta(X_{i_k})+\frac{1}{M}\sum_{k=1}^M \varphi_\theta^\ast(Y_{j_k})$
                \State Update $\theta\in\Theta$ by minimizing $L$ with \texttt{Adam}
            \EndFor
        \EndFor
        \State \textbf{Return} $\nabla\varphi_\theta$
    \end{algorithmic}
\end{algorithm}

\begin{algorithm}
    \caption{Computing the original or sieve convex conjugates}
    \label{alg:convex_conjugate}
    \begin{algorithmic}[1]
        \Require Function $\varphi$; value $y\in\RR^d$; number of epochs $T$; projection radius $M_n$.
        \State Initialize $x\leftarrow 0$
        \State \textbf{def} \texttt{closure}($\varphi$):
            \State \hskip1.5em Compute loss: $l \leftarrow \varphi(x) - \inner{x,y}$
            \State \hskip1.5em \textbf{return} $l$
        \For{$t=1, \cdots, T$}
            \State Update $x$ with $x\leftarrow\texttt{GradientDescent}(\texttt{closure}, x)$
            \State Project $x$ onto $B(0,M_n)$
        \EndFor
        \State Calculate convex conjugate: $\varphi^\ast(y)\leftarrow\inner{x,y}-\varphi(x)$
        \State \textbf{Return} $\varphi^\ast(y)$
    \end{algorithmic}
\end{algorithm}

We summarize our numerical algorithm for the original OT estimator from Equation \eqref{eq: discrete-discrete OT} as in Algorithm~\ref{alg:ot_map_icnn} and Algorithm~\ref{alg:convex_conjugate}. Our end-to-end algorithm is divided into two parts to highlight our method for computing the convex conjugate. For Algorithm~\ref{alg:ot_map_icnn}, we use Adam optimizer \citep{kingma2014adam} due to its superior performance compared to gradient descent in our experiments. For Algorithm~\ref{alg:convex_conjugate}, setting $M_n\leftarrow+\infty$ recovers the original convex conjugate in Equation \eqref{eq: discrete-discrete OT}, while setting $M_n\leftarrow\max_{1\leq i\leq n}\|X_i\|_2$ recovers the sieve convex conjugate in Equation \eqref{eq: new discrete-discrete estimator, overview}.

\subsection{Simulation Results}\label{sec: Simulation Results}

In this section, we demonstrate the performance of OT map estimators from Equation \eqref{eq: discrete-discrete OT} and Equation \eqref{eq: new discrete-discrete estimator, overview}. in the univariate case. Due to space limitations, we refer readers to Appendix A for results in multivariate scenarios ($d=1,2,3,5,10$), along with visualizations in both univariate and bivariate cases.

For the univariate case, we consider two types of distributions for $P$: (1) standard normal distribution, (2) $t$-distribution with 6 degree-of-freedom, $t(6)$. In each case, there are three OT maps to be estimated: (1) the rank function (i.e. the target measure is $Q=U(0,1)$), (2) a linear transformation $\nabla\varphi_0(z)=3z+5$, (3) a signed-quadratic transformation $\nabla\varphi_0(z)=\sign{z}\cdot z^2$.\footnote{$\sign{\cdot}$ is the sign function, i.e. $\sign{x}=1$ for $x>0$, $\sign{0}=0$, and $\sign{x}=-1$ for $x<0$.} The independent empirical measures $\PP_n$ and $\QQ_N$ are randomly sampled for each experiment, with sample sizes set to be $n=N=100,300,500$ and $1000$.

We implement the OT map estimators using \texttt{PyTorch} \citep{paszke2019pytorch}. For the ICNN in Algorithm~\ref{alg:ot_map_icnn}, we adopt the activation function $s(\cdot)$ as Exponential Linear Unit (ELU) \citep{clevert2015fast} for its smoothness. The depth, i.e. the number of layers $K$, is set to be 3. The width, i.e. the number of neurons in the hidden layer, is set to be 15. The learning rate is set to be 0.001 for both optimizers in Algorithm~\ref{alg:ot_map_icnn} and Algorithm~\ref{alg:convex_conjugate}, with both numbers of epochs set to be 500. The mini-batches in Algorithm~\ref{alg:ot_map_icnn} are realized with a \texttt{DataLoader} with shuffling enabled, and $m=M=50$. 

Additionally, we find that the depth, width, and parameter initialization of the ICNN significantly affect performance. While we are unable to guarantee superior performance of the model at this moment, we believe that this engineering issue is beyond the scope of our work. Our goal is not to represent the state-of-the-art for any particular nonparametric model, but to demonstrate the convenience and feasibility of our combination of the objective function from the semi-dual form of the Kantorovich problem, the ICNN model, and the optimization approach for handling the convex conjugate presented in Algorithm~\ref{alg:ot_map_icnn} and Algorithm~\ref{alg:convex_conjugate}. 

To evaluate these OT map estimators, for normal and $t(6)$ distributions, we use additional 1,000 or 100,000 i.i.d. samples to approximate the $L^2(P)$ losses of the estimated OT maps with Monte-Carlo method. In each case, the experiment is repeated independently 100 times. The $L^2(P)$ losses of the original OT estimator, $\nabla\hat\varphi_{n,N}$, from Equation \eqref{eq: discrete-discrete OT} and our new one, $\nabla\tilde\varphi_{n,N}$, from Equation \eqref{eq: new discrete-discrete estimator, overview}, in the univariate case are summarized in Table~\ref{table: numerical results for univariate case}. 

\begin{table}[h]
\scalebox{0.8}{
\begin{tabular}{ccc|cc|cc|cc|cc}
\hline
 & & & \multicolumn{2}{c|}{$n=N=100$} & \multicolumn{2}{c|}{$n=N=300$} & \multicolumn{2}{c|}{$n=N=500$} & \multicolumn{2}{c}{$n=N=1000$} \\ \hline
 Estimators & $P$ & OT map & Mean & SD & Mean & SD & Mean & SD & Mean & SD \\ \hline
\multirow{6}{*}{$\nabla\hat\varphi_{n,N}$}& \multirow{3}{*}{$N(0,1)$} & Rank function & 0.0557 & 0.0245 & 0.0315 & 0.0113 & 0.0242 & 0.0105 & 0.0199 & 0.0075 \\ 
 &  & Linear & 0.7057 & 0.2704 & 0.3590 & 0.1080 & 0.2986 & 0.1045 & 0.2389 & 0.0718 \\ 
 &  & Signed-quadratic & 0.7048 & 0.3519 & 0.4183 & 0.1416 & 0.3478 & 0.1193 & 0.2578 & 0.0872 \\ \cline{2-11}
 & \multirow{3}{*}{$t(6)$} & Rank function & 0.0615 & 0.0317 & 0.0362 & 0.0116 & 0.0299 & 0.0090 & 0.0286 & 0.0069 \\ 
 &  & Linear & 1.5594 & 0.9247 & 0.8131 & 0.3657 & 0.6719 & 0.3122 & 0.4967 & 0.2197 \\ 
 &  & Signed-quadratic & 4.1779 & 3.8534 & 3.0617 & 4.2378 & 2.5647 & 2.6772 & 1.7090 & 0.6826 \\ \hline 
\multirow{6}{*}{$\nabla\tilde\varphi_{n,N}$}& \multirow{3}{*}{$N(0,1)$} & Rank function & 0.0495 & 0.0216 & 0.0311 & 0.0118 & 0.0261 & 0.0102 & 0.0220 & 0.0087 \\ 
 &  & Linear & 0.6816 & 0.2025 & 0.3984 & 0.1298 & 0.3060 & 0.0940 & 0.2192 & 0.0656 \\ 
 &  & Signed-quadratic & 0.6954 & 0.1929 & 0.4232 & 0.1281 & 0.3387 & 0.1051 & 0.2546 & 0.0842 \\ \cline{2-11}
 & \multirow{3}{*}{$t(6)$} & Rank function & 0.0575 & 0.0258 & 0.0366 & 0.0121 & 0.0317 & 0.0098 & 0.0274 & 0.0065 \\ 
 &  & Linear & 1.1227 & 0.2287 & 0.7425 & 0.1427 & 0.5997 & 0.1252 & 0.4679 & 0.0784 \\ 
 &  & Signed-quadratic & 2.4271 & 0.3000 & 2.1068 & 0.4620 & 1.8294 & 0.3308 & 1.6978 & 0.4923 \\ \hline 
\end{tabular}
}
\caption{$L^2(P)$ losses of OT map estimators from Equation \eqref{eq: discrete-discrete OT} and Equation \eqref{eq: new discrete-discrete estimator, overview} in univariate case. As they are obtained via Monte-Carlo method, both the mean and standard deviation (denoted as SD) of the losses are provided.}
\label{table: numerical results for univariate case}
\end{table}

Consistent with our theoretical findings, the $L^2(P)$ losses are generally larger for the $t(6)$ distribution compared to the standard normal distribution. This discrepancy is particularly pronounced when estimating the signed-quadratic OT map, which has a larger parameter $b$ in terms of $(\beta,b)$-smoothness. For the $t(6)$ distribution, this leads to a significant increase in loss, whereas for the normal distribution $P=N(0,1)$, the increase in loss is minimal. This phenomenon aligns with our convergence rates in Theorem \ref{thm: donsker function class with strongly convex, main text} and Theorem \ref{thm: donsker function class without strongly convex, main text}, as a larger parameter $b$ (which is $a$ in Theorem \ref{thm: donsker function class with strongly convex, main text}) results in slower convergence rates for polynomial-tailed distributions. However, this parameter does not affect the convergence rates (up to logarithm terms) for sub-Weibull $P$. 

Furthermore, a distinct performance difference is observed between the two OT map estimators. While increasing the sample size effectively reduces the $L^2(P)$ losses, our new OT estimator $\nabla\tilde\varphi_{n,N}$ consistently outperforms the original estimator $\nabla\hat\varphi_{n,N}$, particularly for the heavy-tailed $t$-distribution. For instance, when estimating the signed-quadratic OT map under $P=t(6)$, the mean loss of $\nabla\tilde\varphi_{n,N}$ with $n=N=100$ is 2.4271, which is comparable to that of $\nabla\hat\varphi_{n,N}$ with $n=N=500$, i.e., 2.5647. Additionally, the new estimator $\nabla\tilde\varphi_{n,N}$ exhibits greater robustness, evidenced by smaller standard deviations of losses, especially in scenarios involving heavy-tailed $t$-distributions.

\section{Conclusion}\label{conclusion}

Optimal transport (OT) map estimation has emerged as a central topic in statistics, applied mathematics, machine learning, and numerous scientific fields. Despite recent advancements, existing methods are often constrained by restrictive regularity assumptions, such as the compactness and convexity of probability measure supports and the bi-Lipschitz property of OT maps. These limitations exclude many practical cases, including fundamental examples like the rank function of a normal or  $t$-distribution and probability measures with non-convex supports.

This work addresses these challenges by bridging the gap between current methods and the theoretical limit established by Brenier's Theorem. We extend OT map estimation to general sub-Weibull and polynomial-tailed distributions without requiring boundedness or convexity of the probability supports. For Donsker function classes, we achieve faster convergence rates by introducing novel Poincaré-type inequalities. Unlike global Poincaré inequalities, which necessitate sub-exponential probability measures, our Poincaré-type inequalities depend only on local density boundedness and mild topological properties of the support. This innovation significantly broadens the range of applicable probability measures. Additionally, our new empirical process results allow us to handle distributions lacking fourth-order moments, pushing OT map estimation closer to the theoretical limit set by Brenier's Theorem.

Moreover, we introduce an alternative plug-in estimator that removes the need for the $(\alpha,a)$-convexity assumption of Brenier's potential. This enables, for the first time, the estimation of the rank function for univariate normal and $t$-distributions using optimal transport theory. Our estimators not only expand the scope of OT map estimation but are also conceptually straightforward and easy to implement. Numerical experiments, leveraging modern neural networks, demonstrate their effectiveness and robustness, particularly in challenging settings involving heavy-tailed distributions.

\bibliographystyle{agsm}
\bibliography{arxiv_main.bib}

\newpage
\appendix

\section*{Appendix A: Numerical Experiments and Visualization}

We first list the algorithm of computing the discreet-discrete OT map estimators (Algorithm~\ref{alg:ot_map_icnn_projected_supp}) with the original convex conjugate (Algorithm~\ref{alg:convex_conjugate_supp}) and the sieved one (Algorithm~\ref{alg:convex_conjugate_projected_supp}) respectively. 

\begin{algorithm}
    \caption{Plug-in OT Map Estimators $\nabla\hat\varphi_{n,N}$ and $\nabla\tilde\varphi_{n,N}$ with ICNN}
    \label{alg:ot_map_icnn_projected_supp}
    \begin{algorithmic}[1]
        \Require Samples $(X_i)_{i=1}^n$ and $(Y_j)_{j=1}^N$; ICNN model $\Fcal_\Theta$; number of epochs $T$; batch sizes $m, M$.
        \State Initialize $\varphi_\theta\in\Fcal_\Theta$
        \State Compute $M_n\leftarrow \max_{i\in[n]}\|X_i\|_2$
        \For{$t=1, \cdots, T$}
            \For{mini-batch $(X_{i_k})_{k=1}^m$ in $(X_i)_{i=1}^n$, and $(Y_{j_k})_{j=1}^M$ in $(Y_j)_{j=1}^N$}
                \State Compute $(\varphi_\theta^\ast(Y_{j_k}))_{k=1}^M$ with Algorithm~\ref{alg:convex_conjugate_projected_supp}
                \State Compute loss: $L\leftarrow\frac{1}{m}\sum_{k=1}^m \varphi_\theta(X_{i_k})+\frac{1}{M}\sum_{k=1}^M \varphi_\theta^\ast(Y_{j_k})$
                \State Update $\theta\in\Theta$ by minimizing $L$ with \texttt{Adam}
            \EndFor
        \EndFor
        \State \textbf{Return} $\nabla\varphi_\theta$
    \end{algorithmic}
\end{algorithm}

\begin{algorithm}
    \caption{Computing the original convex conjugate with gradient descent}
    \label{alg:convex_conjugate_supp}
    \begin{algorithmic}[1]
        \Require Function $\varphi$; value $y\in\RR^d$; number of epochs $T$.
        \State Initialize $x\leftarrow 0$
        \State \textbf{def} \texttt{closure}($\varphi$):
            \State \hskip1.5em Compute loss: $\ell \leftarrow \varphi(x) - \inner{x,y}$
            \State \hskip1.5em \textbf{return} $\ell$

        \State Optimize over $x$: $x\leftarrow\texttt{GradientDescent}(\texttt{closure}, x)$
        \State Calculate convex conjugate: $\varphi^\ast(y)\leftarrow\inner{x,y}-\varphi(x)$
        \State \textbf{Return} $\varphi^\ast(y)$
    \end{algorithmic}
\end{algorithm}

\begin{algorithm}
    \caption{Computing the sieved convex conjugate with projected gradient descent}
    \label{alg:convex_conjugate_projected_supp}
    \begin{algorithmic}[1]
        \Require Function $\varphi$; value $y\in\RR^d$; number of epochs $T$; projection radius $M_n$.
        \State Initialize $x\leftarrow 0$
        \State \textbf{def} \texttt{closure}($\varphi$):
            \State \hskip1.5em Compute loss: $l \leftarrow \varphi(x) - \inner{x,y}$
            \State \hskip1.5em \textbf{return} $l$
        \For{$t=1, \cdots, T$}
            \State Update $x$ with $x\leftarrow\texttt{GradientDescent}(\texttt{closure}, x)$
            \State Project $x$ onto $B(0,M_n)$
        \EndFor
        \State Calculate convex conjugate: $\varphi^\ast(y)\leftarrow\inner{x,y}-\varphi(x)$
        \State \textbf{Return} $\varphi^\ast(y)$
    \end{algorithmic}
\end{algorithm}

As the numerical performance of the OT map estimators in the univariate case was presented in the main paper due to space constraints, we now extend the analysis to different dimensional spaces. Specifically, we consider dimensions $d=1,2,3,5$ and $10$.

For each dimension, we continue to examine two types of distributions for $P$:
\begin{enumerate}
    \item $d$-dimensional standard normal distribution.
    \item $d$-dimensional $t$-distribution with 6 degree-of-freedom.
\end{enumerate}
In these two cases, the coordinates of the random vectors are independently distributed according to either the standard normal distribution or the $t(6)$ distribution.

In the univariate case, we have considered three types of OT maps to estimate: for $z\in\RR$,
\begin{enumerate}
    \item Rank function: The OT is the CDF of either standard normal or $t(6)$ distribution, i.e., the target measure is $Q = U(0,1)$. With some abuse of notation, denote $\nabla\varphi_0(z)=\mathbf{F}(z)$. 
    \item Linear transformation: $\nabla\varphi_0(z) = 3z + 5$.
    \item Signed-quadratic transformation: $\nabla\varphi_0(z) = \text{sign}(z) \cdot z^2$.
\end{enumerate}

For the multivariate case, where the OT maps are functions between two $\RR^d$ spaces, we define the OT maps to be the composition of these three cases: for $\mathbf{z}=(z_1,\cdots,z_d)^\top\in\RR^d$,
\begin{enumerate}
    \item Rank function: 
    $$\nabla\varphi_0(\mathbf{z})=\begin{pmatrix}
        \mathbf{F}(z_1) \\ \mathbf{F}(z_2) \\ \vdots \\ \mathbf{F}(z_d)
    \end{pmatrix}.$$
    \item Linear function: 
    $$\nabla\varphi_0(\mathbf{z})=\begin{pmatrix}
        3z_1 + 5 \\ 3z_2 + 5 \\ \vdots \\ 3z_d + 5
    \end{pmatrix}.$$
    \item Signed-quadratic function:  
    $$\nabla\varphi_0(\mathbf{z})=\begin{pmatrix}
        \text{sign}(z_1) \cdot z_1^2 \\ \text{sign}(z_2) \cdot z_2^2 \\ \vdots \\ \text{sign}(z_d) \cdot z_d^2
    \end{pmatrix}.$$
\end{enumerate}

We refer readers to the main paper for the architecture of the ICNN we adopted and the optimization settings. The independent empirical measures $\PP_n$ and $\QQ_N$ are randomly sampled for each experiment, with sample sizes set to be $n=N=100,300,500$ and $1000$.

To evaluate these OT estimators, for $d$-dimensional normal and $d$-dimensional $t(6)$ distributions, we use additional 1,000 or 100,000 i.i.d. samples to approximate the $L^2(P)$ losses of the estimated OT maps with Monte-Carlo method. In each case, the experiment is repeated independently 100 times. The $L^2(P)$ losses of the original OT estimator, $\nabla\hat\varphi_{n,N}$, from Equation (6) in the main paper and our new one, $\nabla\tilde\varphi_{n,N}$, from Equation (7) in the main paper, are summarized in Table~\ref{table: numerical results of old OT estimator, supp} and Table~\ref{table: numerical results of new OT estimator, supp}, respectively.

\begin{table}[H]
\centering
\scalebox{0.9}{
\begin{tabular}{ccc|cc|cc|cc|cc}
\hline
 & $\nabla\hat\varphi_{n,N}$ & & \multicolumn{2}{c|}{$n=N=100$} & \multicolumn{2}{c|}{$n=N=300$} & \multicolumn{2}{c|}{$n=N=500$} & \multicolumn{2}{c}{$n=N=1000$} \\ \hline
 & $P$ & OT map & Mean & SD & Mean & SD & Mean & SD & Mean & SD \\ \hline
\multirow{6}{*}{$d=1$} & \multirow{3}{*}{$N(0,1)$} & Rank function & 0.0557 & 0.0245 & 0.0315 & 0.0113 & 0.0242 & 0.0105 & 0.0199 & 0.0075 \\ 
 &  & Linear & 0.7057 & 0.2704 & 0.3590 & 0.1080 & 0.2986 & 0.1045 & 0.2389 & 0.0718 \\ 
 &  & Signed-quadratic & 0.7048 & 0.3519 & 0.4183 & 0.1416 & 0.3478 & 0.1193 & 0.2578 & 0.0872 \\ \cline{2-11}
 & \multirow{3}{*}{$t(6)$} & Rank function & 0.0615 & 0.0317 & 0.0362 & 0.0116 & 0.0299 & 0.0090 & 0.0286 & 0.0069 \\ 
 &  & Linear & 1.5594 & 0.9247 & 0.8131 & 0.3657 & 0.6719 & 0.3122 & 0.4967 & 0.2197 \\ 
 &  & Signed-quadratic & 4.1779 & 3.8534 & 3.0617 & 4.2378 & 2.5647 & 2.6772 & 1.7090 & 0.6826 \\ \hline 
\multirow{6}{*}{$d=2$} & \multirow{3}{*}{$N(0,1)$} & Rank function & 0.0692 & 0.0175 & 0.0411 & 0.0089 & 0.0332 & 0.0064 & 0.0253 & 0.0055 \\ 
 &  & Linear & 1.1353 & 0.3010 & 0.5598 & 0.0924 & 0.4379 & 0.0718 & 0.3292 & 0.0486 \\ 
 &  & Signed-quadratic & 0.7336 & 0.1983 & 0.4760 & 0.1214 & 0.3849 & 0.0892 & 0.3027 & 0.0798 \\ \cline{2-11}
 & \multirow{3}{*}{$t(6)$} & Rank function & 0.0728 & 0.0198 & 0.0444 & 0.0083 & 0.0356 & 0.0057 & 0.0287 & 0.0046 \\ 
 &  & Linear & 2.1658 & 0.9524 & 1.1784 & 0.4359 & 0.8302 & 0.1969 & 0.6366 & 0.1450 \\ 
 &  & Signed-quadratic & 3.9665 & 3.0174 & 3.2203 & 1.7875 & 2.6183 & 1.2406 & 2.1947 & 1.5544 \\ \hline 
\multirow{6}{*}{$d=3$} & \multirow{3}{*}{$N(0,1)$} & Rank function & 0.0799 & 0.0131 & 0.0473 & 0.0064 & 0.0390 & 0.0059 & 0.0301 & 0.0045 \\ 
 &  & Linear & 1.5343 & 0.3312 & 0.8131 & 0.0997 & 0.6177 & 0.0779 & 0.4441 & 0.0478 \\ 
 &  & Signed-quadratic & 0.7355 & 0.1683 & 0.5241 & 0.0920 & 0.4445 & 0.0792 & 0.3610 & 0.0482 \\ \cline{2-11}
 & \multirow{3}{*}{$t(6)$} & Rank function & 0.0856 & 0.0140 & 0.0525 & 0.0072 & 0.0406 & 0.0047 & 0.0322 & 0.0032 \\ 
 &  & Linear & 2.3592 & 0.6550 & 1.5799 & 0.4178 & 1.1217 & 0.2206 & 0.7740 & 0.1080 \\ 
 &  & Signed-quadratic & 3.5714 & 1.8129 & 3.7205 & 5.0072 & 2.9627 & 2.4499 & 2.1041 & 1.1064 \\ \hline 
\multirow{6}{*}{$d=5$} & \multirow{3}{*}{$N(0,1)$} & Rank function & 0.0958 & 0.0099 & 0.0658 & 0.0053 & 0.0521 & 0.0046 & 0.0387 & 0.0032 \\ 
 &  & Linear & 2.2026 & 0.2667 & 1.2122 & 0.0986 & 0.9622 & 0.0797 & 0.6562 & 0.0451 \\ 
 &  & Signed-quadratic & 0.8193 & 0.1365 & 0.6677 & 0.0683 & 0.5988 & 0.0440 & 0.5249 & 0.0328 \\ \cline{2-11}
 & \multirow{3}{*}{$t(6)$} & Rank function & 0.1065 & 0.0100 & 0.0680 & 0.0061 & 0.0522 & 0.0039 & 0.0410 & 0.0028 \\ 
 &  & Linear & 2.7752 & 0.3688 & 2.0871 & 0.3771 & 1.4250 & 0.2199 & 1.0067 & 0.1370 \\ 
 &  & Signed-quadratic & 3.6442 & 1.5187 & 3.4569 & 2.0479 & 3.5799 & 2.7394 & 2.4031 & 1.5449 \\ \hline 
\multirow{6}{*}{$d=10$} & \multirow{3}{*}{$N(0,1)$} & Rank function & 0.1330 & 0.0078 & 0.0998 & 0.0033 & 0.0885 & 0.0034 & 0.0636 & 0.0048 \\ 
 &  & Linear & 2.7187 & 0.1822 & 2.0264 & 0.1074 & 1.5856 & 0.0643 & 1.1139 & 0.0333 \\ 
 &  & Signed-quadratic & 1.0868 & 0.0963 & 0.9777 & 0.0665 & 0.8965 & 0.0481 & 0.8085 & 0.0298 \\ \cline{2-11}
 & \multirow{3}{*}{$t(6)$} & Rank function & 0.1431 & 0.0069 & 0.1110 & 0.0050 & 0.0884 & 0.0065 & 0.0559 & 0.0022 \\ 
 &  & Linear & 2.9331 & 0.2062 & 2.9798 & 0.2964 & 2.2717 & 0.1918 & 1.5849 & 0.0993 \\ 
 &  & Signed-quadratic & 4.5246 & 2.1171 & 4.5407 & 2.0494 & 3.9708 & 1.5751 & 3.0236 & 0.9330 \\ \hline 
\end{tabular}
}
\caption{$L^2(P)$ losses of the OT estimator $\nabla\hat\varphi_{n,N}$ from Equation (6) in the main paper. As they are obtained via Monte-Carlo method, both the mean and standard deviation (denoted as SD) of the losses are provided.}
\label{table: numerical results of old OT estimator, supp}
\end{table}

\begin{table}[H]
\centering
\scalebox{0.9}{
\begin{tabular}{ccc|cc|cc|cc|cc}
\hline
 & $\nabla\tilde\varphi_{n,N}$ & & \multicolumn{2}{c|}{$n=N=100$} & \multicolumn{2}{c|}{$n=N=300$} & \multicolumn{2}{c|}{$n=N=500$} & \multicolumn{2}{c}{$n=N=1000$} \\ \hline
 & $P$ & OT map & Mean & SD & Mean & SD & Mean & SD & Mean & SD \\ \hline
\multirow{6}{*}{$d=1$} & \multirow{3}{*}{$N(0,1)$} & Rank function & 0.0495 & 0.0216 & 0.0311 & 0.0118 & 0.0261 & 0.0102 & 0.0220 & 0.0087 \\ 
 &  & Linear & 0.6816 & 0.2025 & 0.3984 & 0.1298 & 0.3060 & 0.0940 & 0.2192 & 0.0656 \\ 
 &  & Signed-quadratic & 0.6954 & 0.1929 & 0.4232 & 0.1281 & 0.3387 & 0.1051 & 0.2546 & 0.0842 \\ \cline{2-11}
 & \multirow{3}{*}{$t(6)$} & Rank function & 0.0575 & 0.0258 & 0.0366 & 0.0121 & 0.0317 & 0.0098 & 0.0274 & 0.0065 \\ 
 &  & Linear & 1.1227 & 0.2287 & 0.7425 & 0.1427 & 0.5997 & 0.1252 & 0.4679 & 0.0784 \\ 
 &  & Signed-quadratic & 2.4271 & 0.3000 & 2.1068 & 0.4620 & 1.8294 & 0.3308 & 1.6978 & 0.4923 \\ \hline 
\multirow{6}{*}{$d=2$} & \multirow{3}{*}{$N(0,1)$} & Rank function & 0.0645 & 0.0154 & 0.0400 & 0.0091 & 0.0307 & 0.0058 & 0.0249 & 0.0048 \\ 
 &  & Linear & 0.8107 & 0.1371 & 0.4828 & 0.0872 & 0.3917 & 0.0624 & 0.3011 & 0.0416 \\ 
 &  & Signed-quadratic & 0.6611 & 0.1354 & 0.4590 & 0.1047 & 0.3680 & 0.0845 & 0.2834 & 0.0549 \\ \cline{2-11}
 & \multirow{3}{*}{$t(6)$} & Rank function & 0.0724 & 0.0167 & 0.0436 & 0.0094 & 0.0348 & 0.0064 & 0.0296 & 0.0046 \\ 
 &  & Linear & 1.2575 & 0.1960 & 0.8054 & 0.1090 & 0.6738 & 0.0800 & 0.5383 & 0.0632 \\ 
 &  & Signed-quadratic & 2.4201 & 0.4221 & 1.9977 & 0.2789 & 1.8119 & 0.2981 & 1.7630 & 0.9320 \\ \hline 
\multirow{6}{*}{$d=3$} & \multirow{3}{*}{$N(0,1)$} & Rank function & 0.0772 & 0.0123 & 0.0462 & 0.0063 & 0.0378 & 0.0054 & 0.0291 & 0.0040 \\ 
 &  & Linear & 0.9835 & 0.1285 & 0.5875 & 0.0742 & 0.4990 & 0.0630 & 0.3758 & 0.0372 \\ 
 &  & Signed-quadratic & 0.7057 & 0.1325 & 0.4980 & 0.0668 & 0.4378 & 0.0651 & 0.3635 & 0.0419 \\ \cline{2-11}
 & \multirow{3}{*}{$t(6)$} & Rank function & 0.0823 & 0.0117 & 0.0488 & 0.0056 & 0.0415 & 0.0053 & 0.0329 & 0.0036 \\ 
 &  & Linear & 1.4409 & 0.1969 & 0.8861 & 0.0897 & 0.7731 & 0.1025 & 0.6251 & 0.0676 \\ 
 &  & Signed-quadratic & 2.4398 & 0.4205 & 2.0213 & 0.3743 & 1.8480 & 0.2582 & 1.6798 & 0.3702 \\ \hline 
\multirow{6}{*}{$d=5$} & \multirow{3}{*}{$N(0,1)$} & Rank function & 0.0934 & 0.0085 & 0.0619 & 0.0050 & 0.0499 & 0.0039 & 0.0386 & 0.0029 \\ 
 &  & Linear & 1.2597 & 0.1021 & 0.8297 & 0.0582 & 0.6896 & 0.0412 & 0.5553 & 0.0359 \\ 
 &  & Signed-quadratic & 0.7871 & 0.0861 & 0.6306 & 0.0572 & 0.5805 & 0.0398 & 0.5154 & 0.0278 \\ \cline{2-11}
 & \multirow{3}{*}{$t(6)$} & Rank function & 0.1026 & 0.0079 & 0.0634 & 0.0058 & 0.0514 & 0.0035 & 0.0405 & 0.0025 \\ 
 &  & Linear & 1.7439 & 0.1583 & 1.1517 & 0.0925 & 0.9803 & 0.0717 & 0.8328 & 0.0553 \\ 
 &  & Signed-quadratic & 2.5232 & 0.4057 & 2.1447 & 0.3413 & 2.0146 & 0.3773 & 1.8229 & 0.3984 \\ \hline 
\multirow{6}{*}{$d=10$} & \multirow{3}{*}{$N(0,1)$} & Rank function & 0.1249 & 0.0064 & 0.0964 & 0.0036 & 0.0844 & 0.0047 & 0.0604 & 0.0039 \\ 
 &  & Linear & 1.7548 & 0.1147 & 1.3402 & 0.0530 & 1.1445 & 0.0327 & 0.9283 & 0.0315 \\ 
 &  & Signed-quadratic & 1.0101 & 0.0774 & 0.8678 & 0.0556 & 0.8162 & 0.0417 & 0.7554 & 0.0273 \\ \cline{2-11}
 & \multirow{3}{*}{$t(6)$} & Rank function & 0.1384 & 0.0053 & 0.1051 & 0.0048 & 0.0819 & 0.0059 & 0.0557 & 0.0018 \\ 
 &  & Linear & 2.2381 & 0.1527 & 1.8178 & 0.0853 & 1.5445 & 0.0551 & 1.3009 & 0.0460 \\ 
 &  & Signed-quadratic & 2.8338 & 0.2913 & 2.5350 & 0.3327 & 2.3599 & 0.2395 & 2.1654 & 0.2265 \\ \hline 
\end{tabular}
}
\caption{$L^2(P)$ losses of the OT estimator $\nabla\tilde\varphi_{n,N}$ from Equation (7) in the main paper. As they are obtained via Monte-Carlo method, both the mean and standard deviation (denoted as SD) of the losses are provided.}
\label{table: numerical results of new OT estimator, supp}
\end{table}

In analyzing the results from these two tables, several key observations can be made:
\begin{enumerate}
    \item \textbf{Effect of Sample Size:} Increasing the sample size consistently decreases the $L^2(P)$ losses and enhances the robustness of the OT map estimators, as evidenced by the reduction in standard deviation. 
    
    \item \textbf{Impact of Dimensionality:} As the dimension $d$ increases, the losses tend to increase as well, which is expected due to the ``curse of dimensionality''.
    
    \item \textbf{Comparison of Distributions:} The $L^2(P)$ losses are generally larger for the $t(6)$ distribution than for the normal distribution. This discrepancy is particularly pronounced when estimating the signed-quadratic OT map, which has a larger $b$ parameter in terms of $(\beta,b)$-smoothness. For the $t(6)$ distribution, this leads to a significant increase in loss, whereas for the normal distribution $P=N(0,1)$, the increase in loss is minimal. This observation aligns with the theoretical convergence rates discussed in the main paper, where a larger $b$ parameter results in slower convergence rates for polynomial-tail distributions, while the convergence rates for sub-Weibull distributions remain unaffected (up to logarithmic terms).

    \item \textbf{Comparison of Estimators:} When comparing the two estimators, $\nabla\tilde\varphi_{n,N}$ and $\nabla\hat\varphi_{n,N}$, it is evident that the new estimator $\nabla\tilde\varphi_{n,N}$ is more efficient and robust. For example, in the univariate case for the signed-quadratic OT map under $P = t(6)$, the mean loss for $\nabla\tilde\varphi_{n,N}$ with $n=N=100$ is 2.4271, which is comparable to the mean loss, 2.5647, for $\nabla\hat\varphi_{n,N}$ with $n=N=500$. Moreover, the new estimator not only reduces the mean loss but also exhibits greater robustness, as demonstrated by the smaller standard deviation of its losses, particularly in scenarios involving heavy-tailed $t$-distributions.
\end{enumerate}

To gain a clearer understanding of the performance of the OT map estimators proposed in our study, we have also visualized the estimated OT maps for both the univariate $(d=1)$ and bivariate $(d=2)$ cases. Importantly, to avoid any potential accusations of cherry-picking, the visualizations are based on the first trained model across all settings.

\begin{figure}[H]
  \centering
  \includegraphics[trim=0cm 7cm 0cm 12cm, width=1.4\textwidth, center]{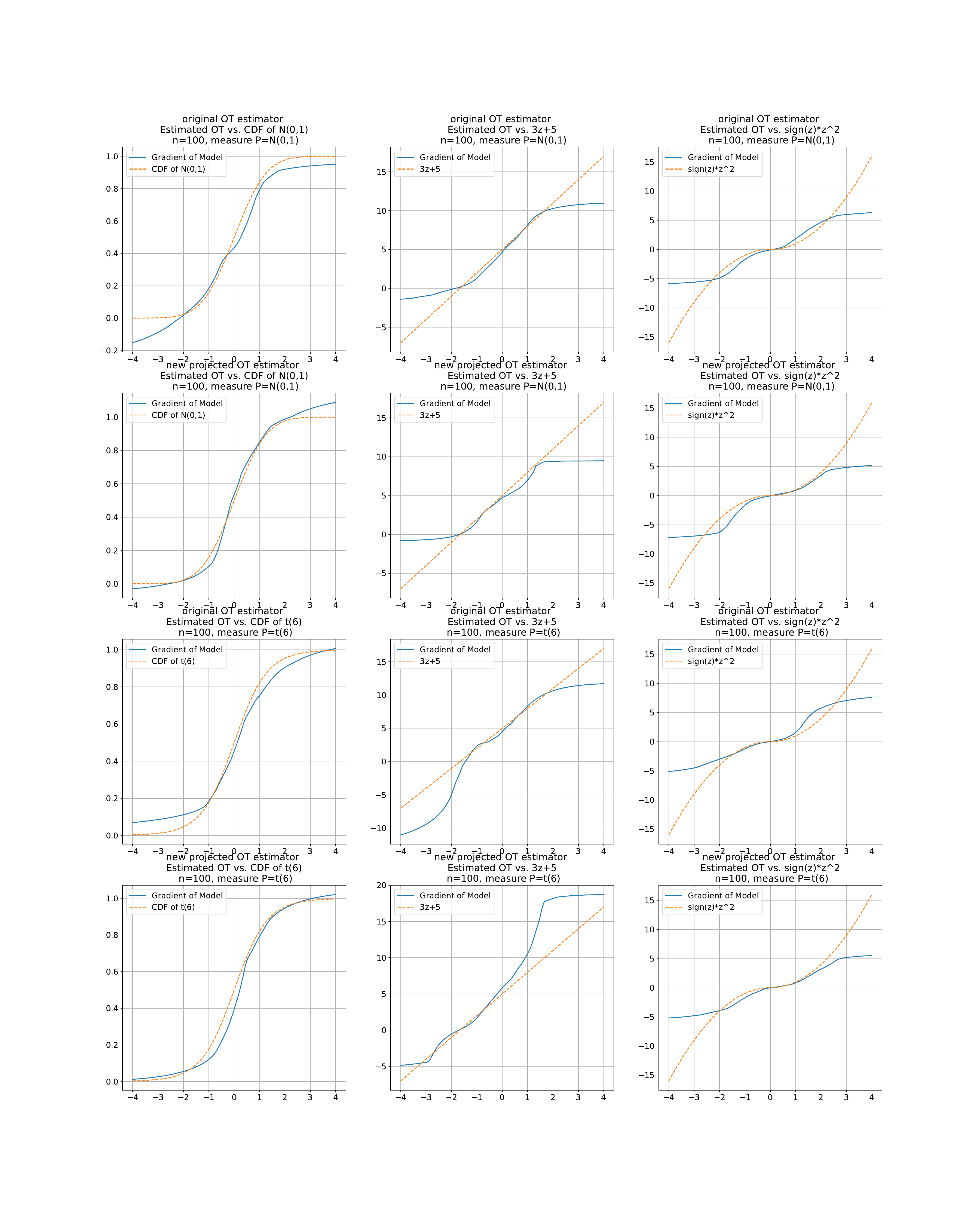}
  \caption{Visualization of univariate estimated OT maps in various settings for sample size $n=N=100$. Best view in color. In each plot, the orange dashed line is the OT map to be estimated and the blue solid line is the estimated OT map as the gradient of an ICNN.}
  \label{fig:n=100}
\end{figure}

\begin{figure}[H]
  \centering
  \includegraphics[trim=0cm 7cm 0cm 12cm, width=1.4\textwidth, center]{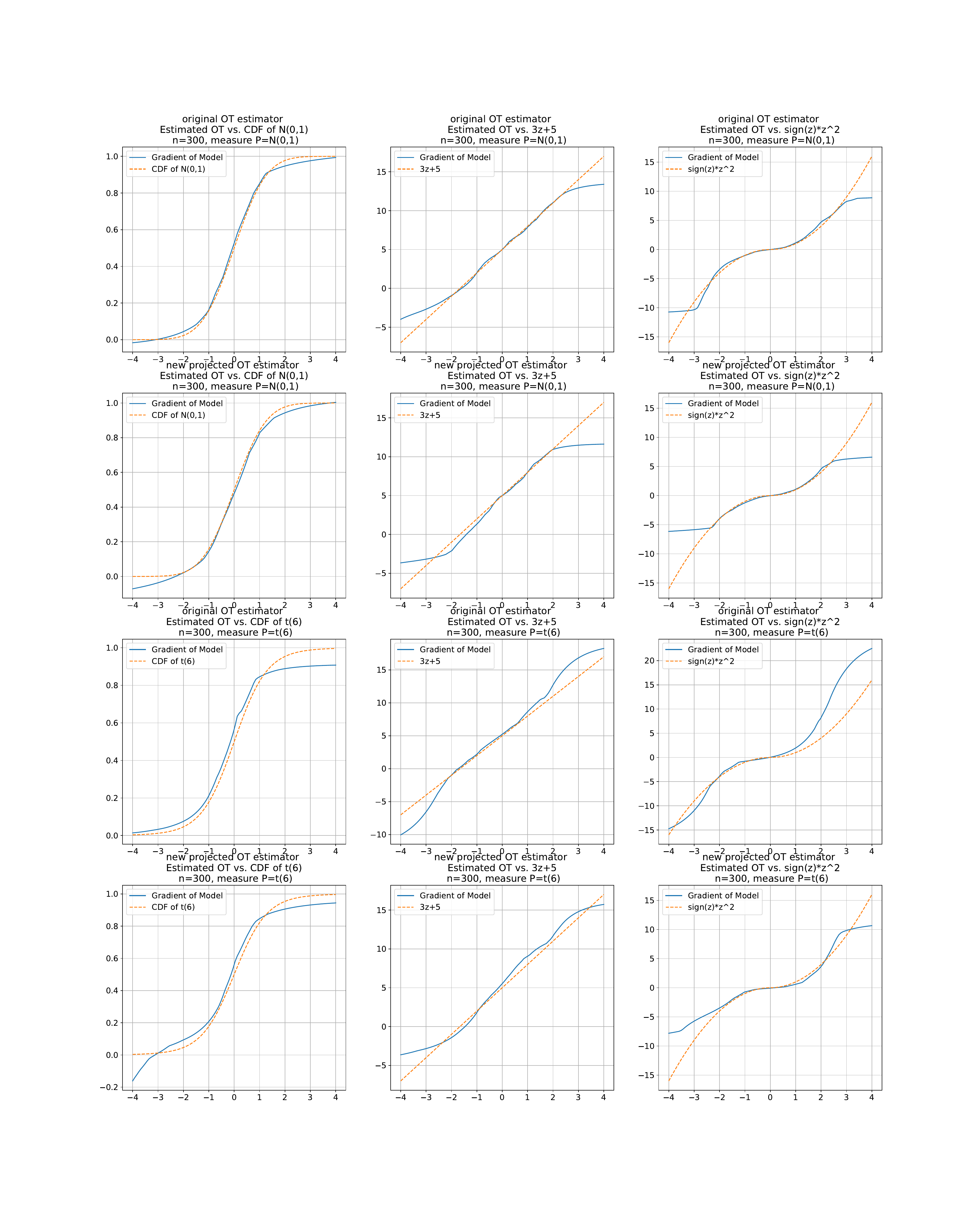}
  \caption{Visualization of univariate estimated OT maps in various settings for sample size $n=N=300$. Best view in color. In each plot, the orange dashed line is the OT map to be estimated and the blue solid line is the estimated OT map as the gradient of an ICNN.}
  \label{fig:n=300}
\end{figure}

\begin{figure}[H]
  \centering
  \includegraphics[trim=0cm 7cm 0cm 12cm, width=1.4\textwidth, center]{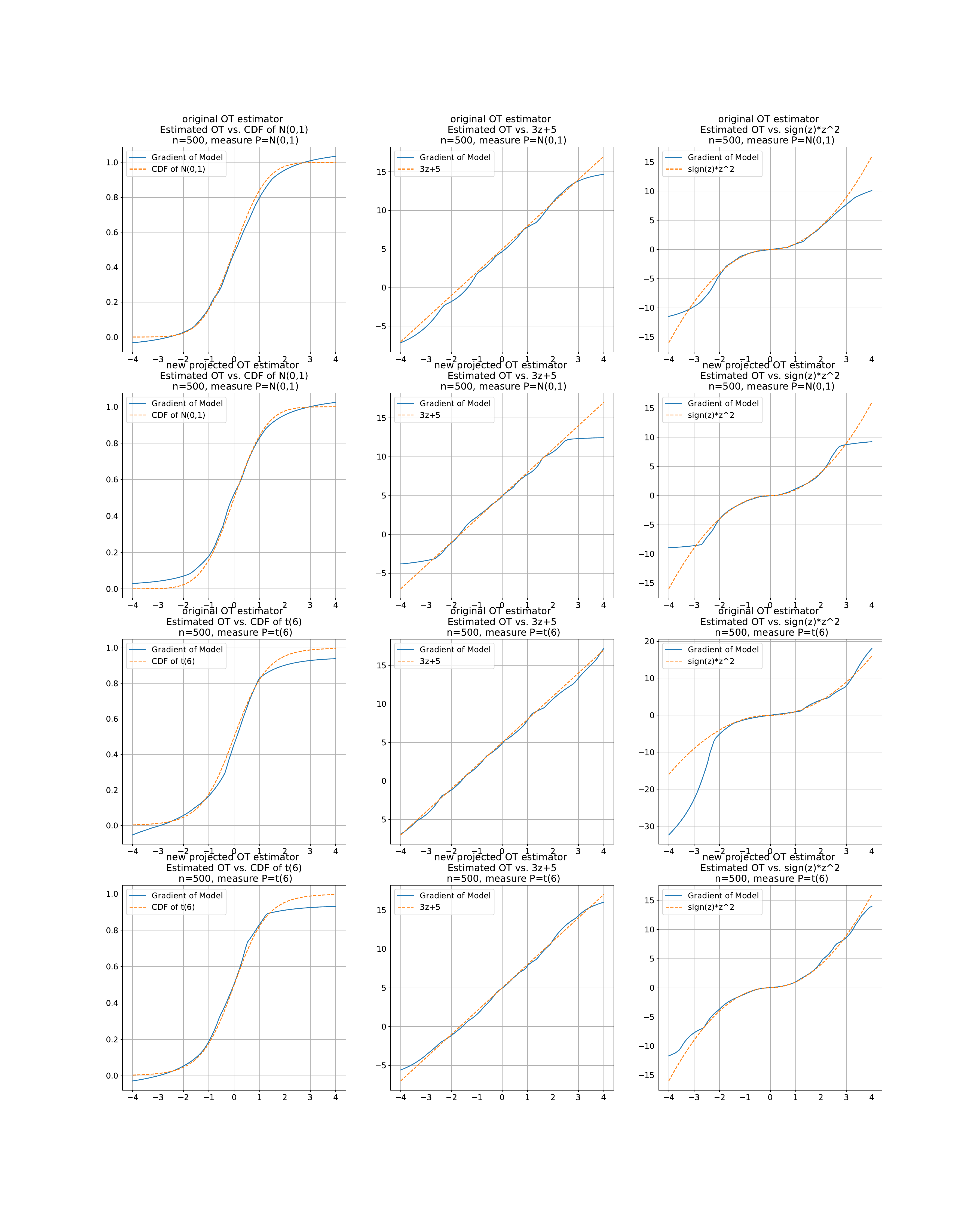}
  \caption{Visualization of univariate estimated OT maps in various settings for sample size $n=N=500$. Best view in color. In each plot, the orange dashed line is the OT map to be estimated and the blue solid line is the estimated OT map as the gradient of an ICNN.}
  \label{fig:n=500}
\end{figure}

\begin{figure}[H]
  \centering
  \includegraphics[trim=0cm 7cm 0cm 12cm, width=1.4\textwidth, center]{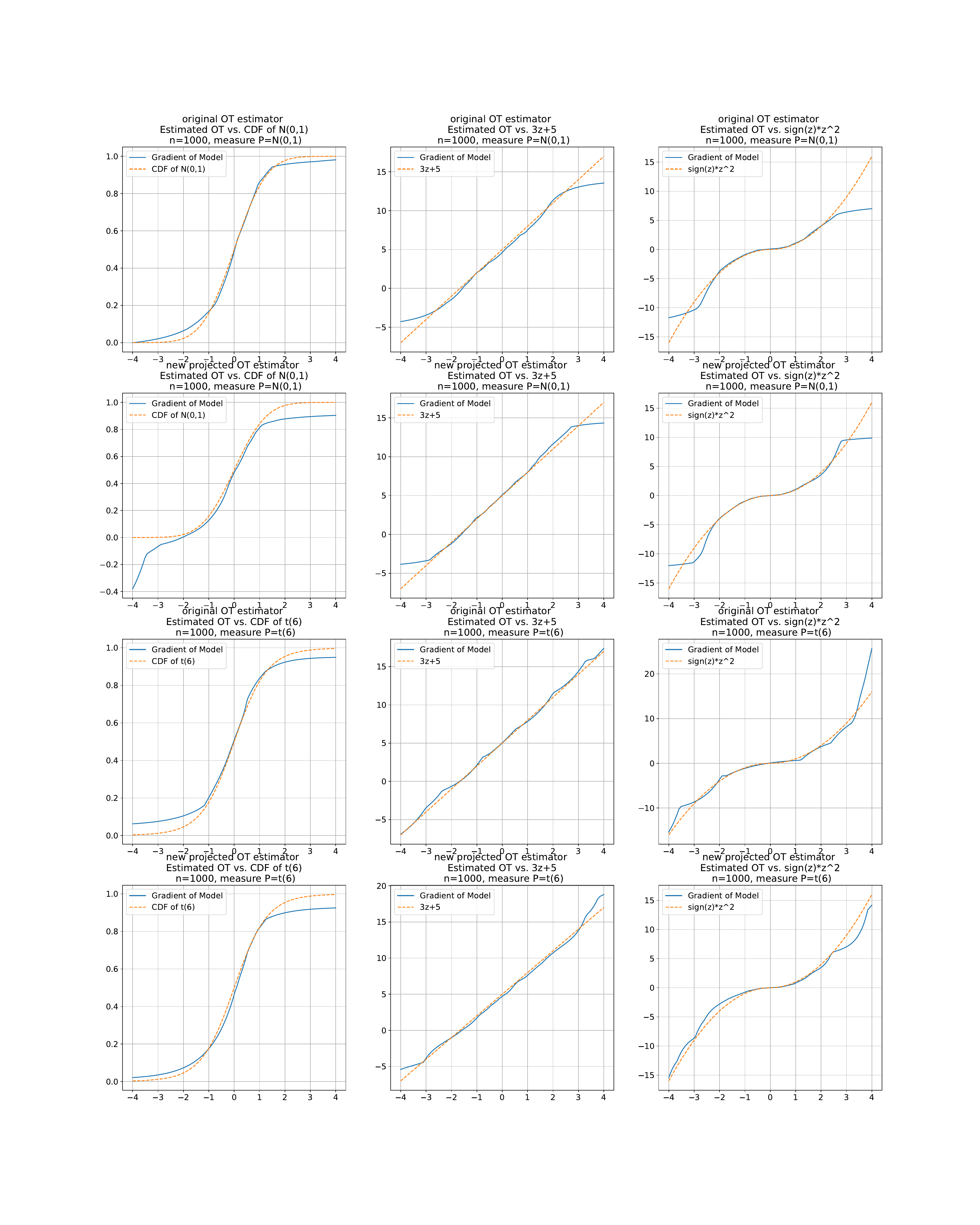}
  \caption{Visualization of univariate estimated OT maps in various settings for sample size $n=N=1000$. Best view in color. In each plot, the orange dashed line is the OT map to be estimated and the blue solid line is the estimated OT map as the gradient of an ICNN.}
  \label{fig:n=1000}
\end{figure}

From Figure~\ref{fig:n=100} to Figure~\ref{fig:n=1000}, we present the univariate estimated OT maps as blue solid lines alongside the ground truth OT maps as orange dashed lines. The following observations can be made from these visualizations:
\begin{enumerate}
    \item As the sample size increases from $n=N=100$ to $n=N=1000$, the estimated OT maps progressively fit the ground truth more accurately. This improvement is consistent with our numerical results and theoretical conclusions, which show a decrease in $L^2(P)$ losses as sample size increases, leading to better performance of the estimators.
    \item The OT map estimators exhibit good performance in regions where the distribution $P$ has a large probability density. In the outer regions, where $P$ has a smaller probability mass, the OT map estimators tend to perform poorly. Specifically, the estimated OT maps $\nabla\tilde\varphi_{n,N}(x)$ and $\nabla\hat\varphi_{n,N}(x)$ approach some constants as $x$ tends to infinity. This behavior is hypothesized to be related to the structure of the ICNN used in our implementation. The Exponential Linear Unit (ELU) activation function, which approximates the ReLU function, tends to behave linearly as $x$ increases towards $+\infty$. Consequently, the ICNN may impose a form of regularization that results in the estimated OT maps (as the gradient of these ICNNs) flattening out in these extrapolation regions. Given this behavior, it is advisable to consider alternative network architectures or non-parametric models according to the prior knowledge of the true OT map.
\end{enumerate}

\begin{figure}[H]
  \centering
  \includegraphics[trim=0cm 7cm 0cm 12cm, width=1.4\textwidth, center]{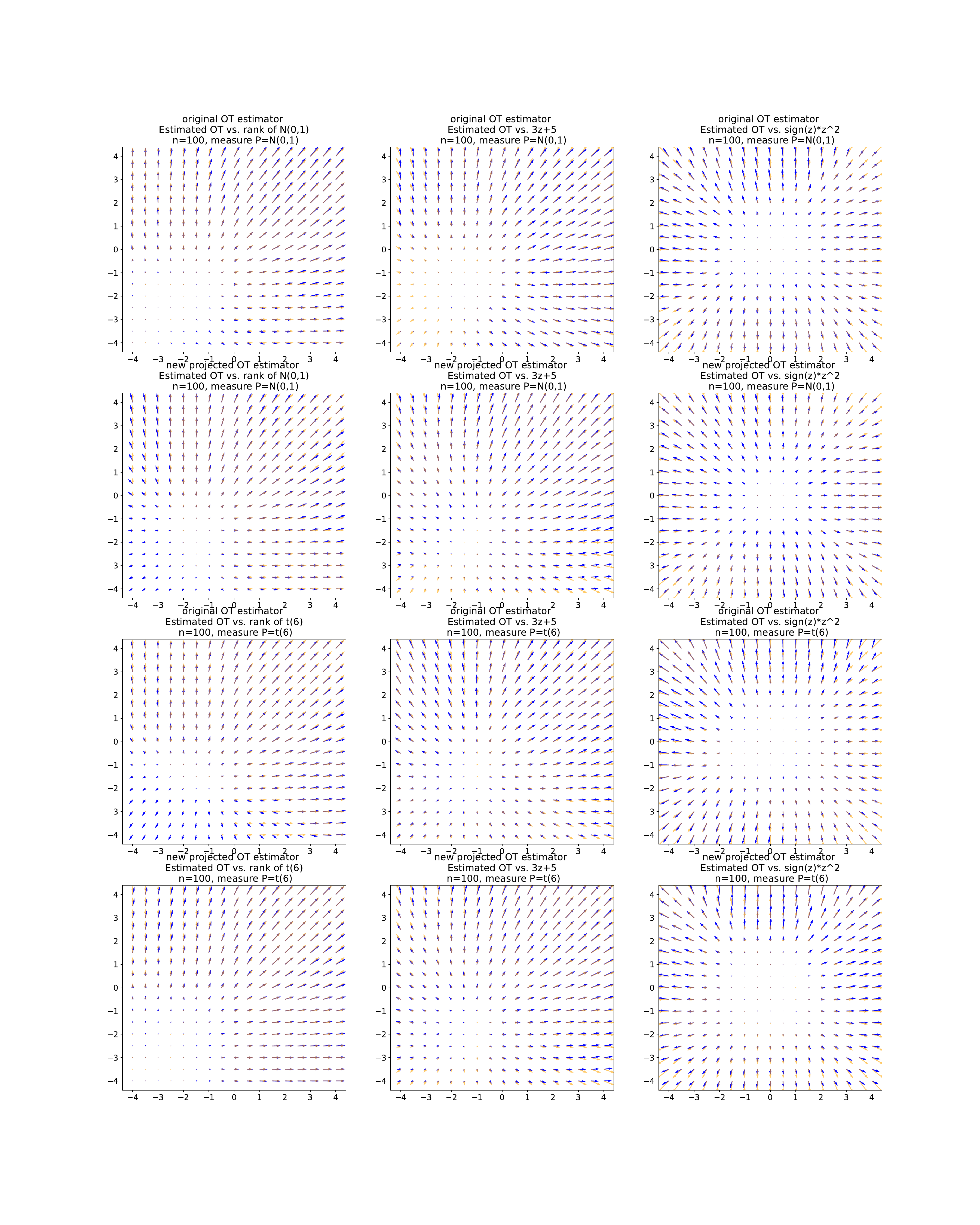}
  \caption{Visualization of 2-dimensional estimated OT maps in various settings for sample size $n=N=100$. Best view in color. In each plot, the orange vector field is the OT map to be estimated and the blue vector field is the estimated OT map as the gradient of an ICNN.}
  \label{fig:d=2_n=100}
\end{figure}

\begin{figure}[H]
  \centering
  \includegraphics[trim=0cm 7cm 0cm 12cm, width=1.4\textwidth, center]{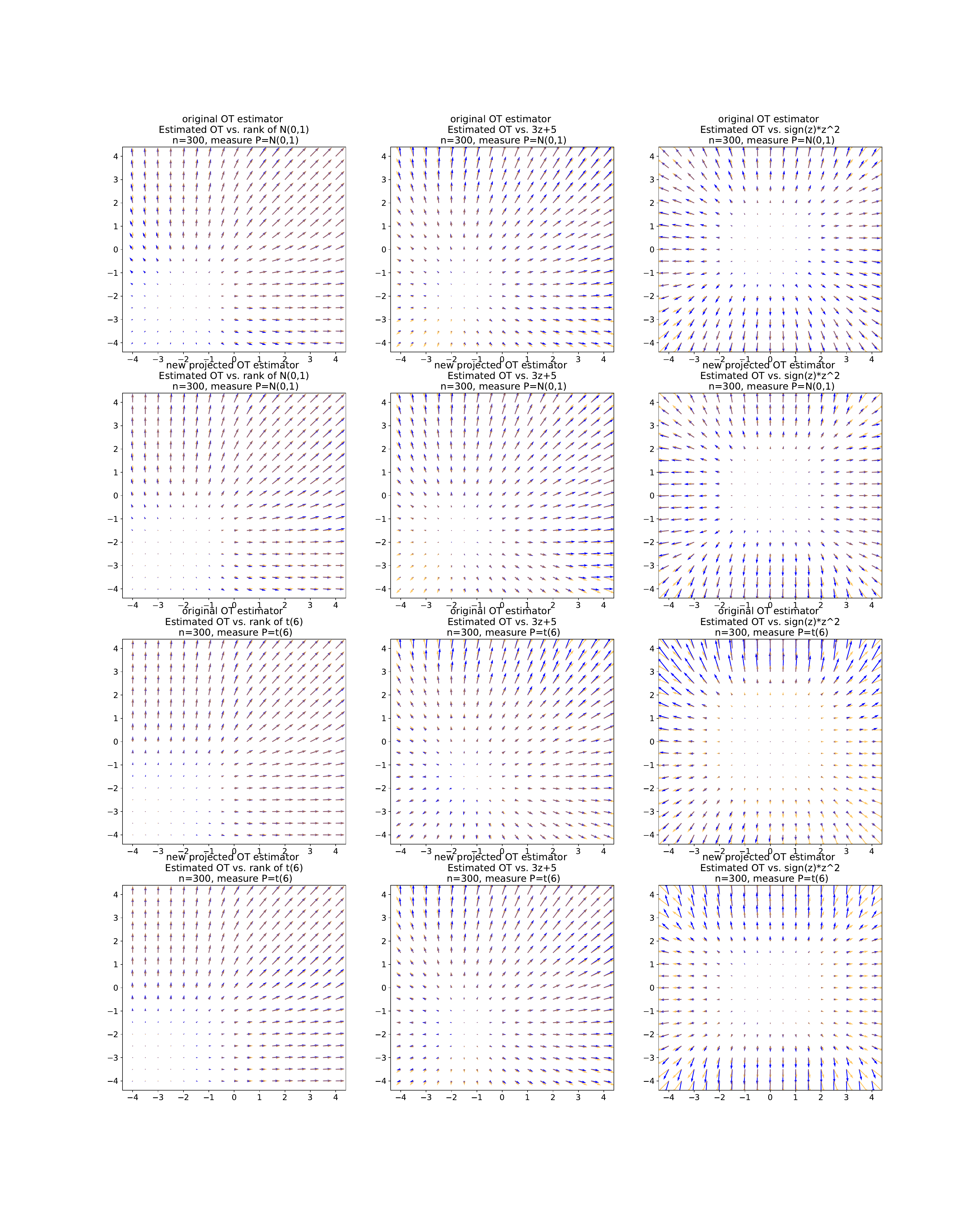}
  \caption{Visualization of 2-dimensional estimated OT maps in various settings for sample size $n=N=300$. Best view in color. In each plot, the orange vector field is the OT map to be estimated and the blue vector field is the estimated OT map as the gradient of an ICNN.}
  \label{fig:d=2_n=300}
\end{figure}

\begin{figure}[H]
  \centering
  \includegraphics[trim=0cm 7cm 0cm 12cm, width=1.4\textwidth, center]{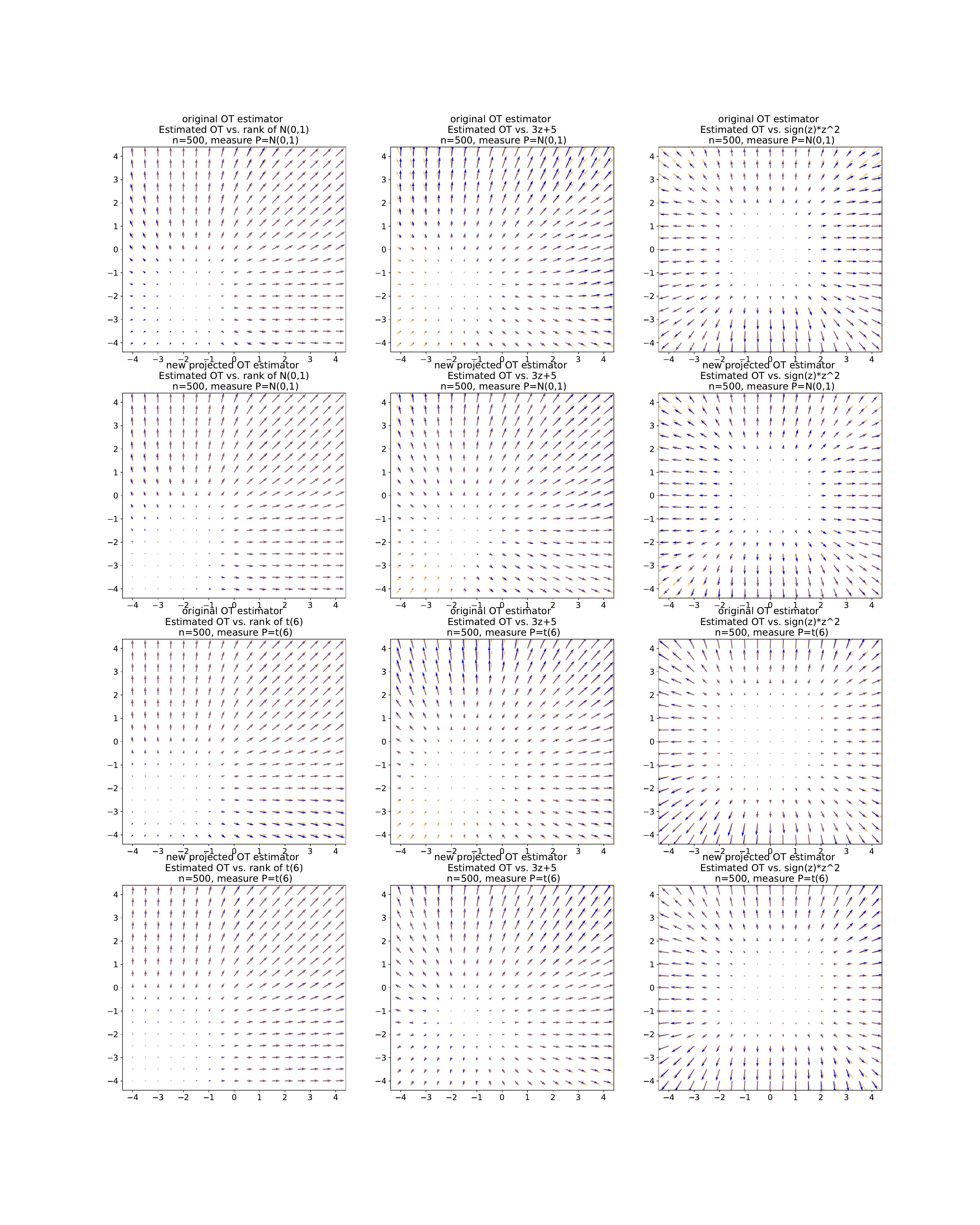}
  \caption{Visualization of 2-dimensional estimated OT maps in various settings for sample size $n=N=500$. Best view in color. In each plot, the orange vector field is the OT map to be estimated and the blue vector field is the estimated OT map as the gradient of an ICNN.}
  \label{fig:d=2_n=500}
\end{figure}

\begin{figure}[H]
  \centering
  \includegraphics[trim=0cm 7cm 0cm 12cm, width=1.4\textwidth, center]{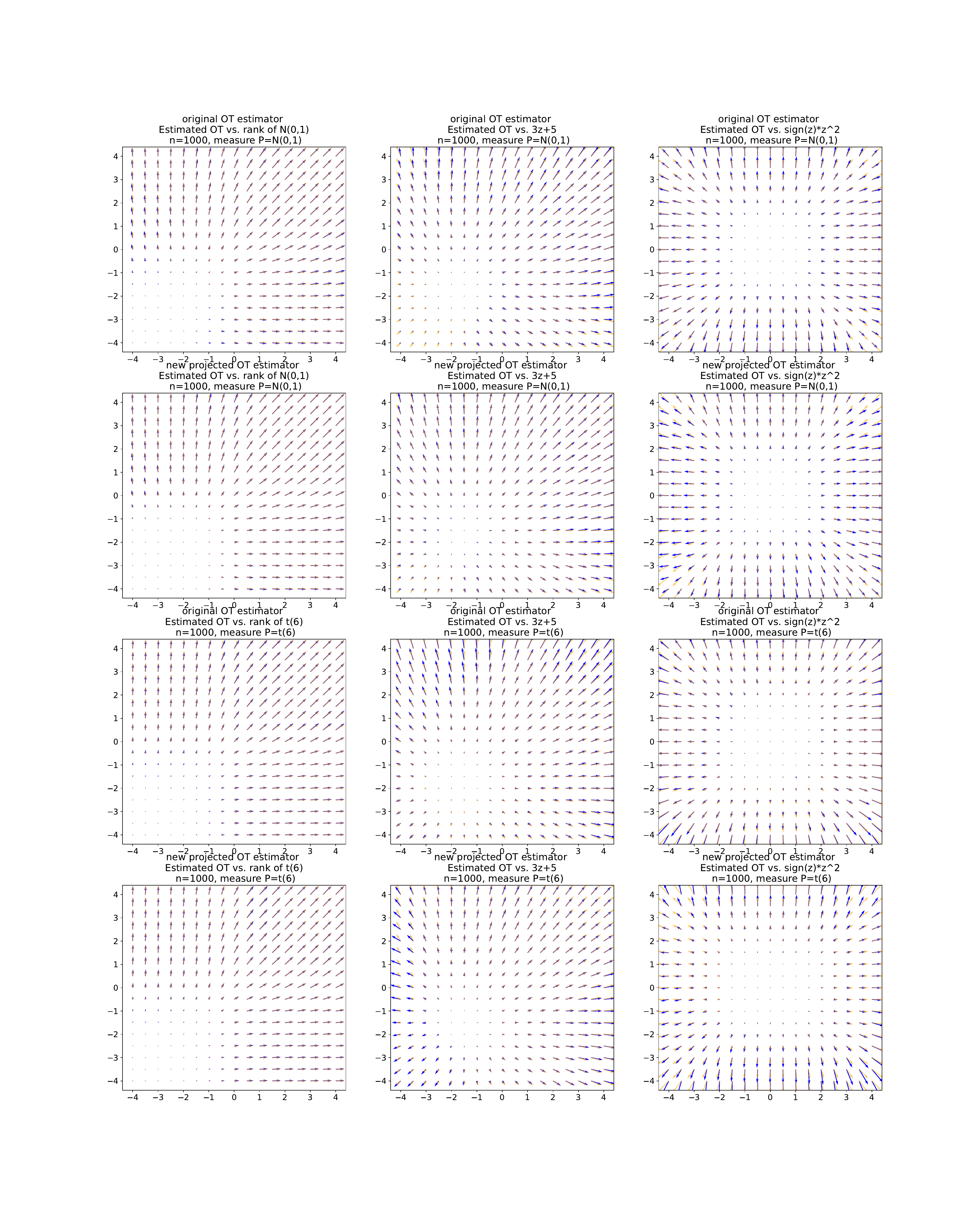}
  \caption{Visualization of 2-dimensional estimated OT maps in various settings for sample size $n=N=1000$. Best view in color. In each plot, the orange vector field is the OT map to be estimated and the blue vector field is the estimated OT map as the gradient of an ICNN.}
  \label{fig:d=2_n=1000}
\end{figure}

In Figure~\ref{fig:d=2_n=100} through Figure~\ref{fig:d=2_n=1000}, we visualize the estimated bivariate OT maps in the form of vector fields, since the OT map is the gradient of a convex function according to Brenier's Theorem.

\end{document}